\documentclass[a4paper,fleqn]{cas-sc}

\usepackage[numbers]{natbib}

\usepackage{amssymb,amsmath,color,amsfonts}
\usepackage{graphicx}
\usepackage{caption}
\usepackage{graphicx} 
\usepackage{array}

\usepackage{mathtools}

\usepackage[shortlabels]{enumitem}

\newtheorem{assumption}{Assumption}
\newtheorem{definition}{Definition}

\newtheorem{lemma}{Lemma}
\newtheorem{theorem}{Theorem}

\newtheorem{proposition}{Proposition}
\newtheorem{remark}{Remark}

\newtheorem{example}{Example}
\newtheorem{problem}{Problem}

\newcommand{\blue}[1]{{\color{black} #1}}
\newcommand{\red}[1]{{\color{blue} #1}}

\newcommand{\R}{\mathbb{R}}

\def\tsc#1{\csdef{#1}{\textsc{\lowercase{#1}}\xspace}}
\tsc{WGM}
\tsc{QE}
\tsc{EP}
\tsc{PMS}
\tsc{BEC}
\tsc{DE}

\newcommand{\mc}{\mathcal}

\newcommand{\B}{\mathbb{B}} 					

\newcommand{\real}{\mathbb{R}} 
\newcommand{\integ}{\mathbb{Z}}

\newcommand{\realpos}{\mathbb{R}_{\geq 0}}

\newcommand{\tsp}{\mathsf{T}} 
 
\newcommand{\inv}{{\negat 1}} 
\newcommand{\negat}{\scalebox{0.75}[.9]{\( - \)}}

\newcommand\oprocendsymbol{\hbox{$\square$}}
\newcommand\oprocend{\relax\ifmmode\else\unskip\hfill%
\fi\oprocendsymbol}

\newcommand{\map}[3]{#1: #2 \rightarrow #3}

\newcommand{\norm}[1]{\left\lVert#1\right\rVert}

\begin{document}
\let\WriteBookmarks\relax
\def\floatpagepagefraction{1}
\def\textpagefraction{.001}
\shorttitle{Online Optimization of LTI Systems Under Persistent Attacks}
\shortauthors{Galarza-Jimenez et~al.}
\title [mode = title]{Online Optimization of LTI Systems Under Persistent Attacks: Stability, Tracking, and Robustness}

\author{Felipe Galarza-Jimenez}

\ead{felipe.galarzajimenez@colorado.edu}

\address{Department of Electrical, Computer, and Energy Engineering, University of Colorado, Boulder, CO, USA.}

\author{Gianluca Bianchin}

\ead{gianluca.bianchin@colorado.edu}

\author{Jorge I. Poveda}

\ead{jorge.poveda@colorado.edu}

\author{Emiliano Dall'Anese}
\ead{emiliano.dallanese@colorado.edu}

\nonumnote{This work was supported by the National Science Foundation (NSF) through the Awards CMMI 2044946 and CNS 1947613, and by the National Renewable Energy Laboratory through the contract UGA-0-41026-148.
  }

\begin{abstract}
We study the stability properties \blue{of a control system composed of a dynamical plant and a feedback controller, the latter generating control signals that can be compromised by a malicious attacker.} We consider two classes of feedback controllers: a static output-feedback controller, and a dynamical gradient-flow controller that seeks to steer the output of the plant towards the solution of a convex optimization problem. In both cases, we analyze the stability properties of the closed-loop system under a class of switching attacks that persistently modify the control inputs generated by the controllers. Our stability analysis leverages the framework of hybrid dynamical systems, Lyapunov-based arguments for switching systems with unstable modes, and singular perturbation theory. Our results reveal that, under a suitable time-scale separation between plant and controllers, the stability of the interconnected system can be preserved when the attack occurs with ``sufficiently low frequency'' in any bounded time interval. We present simulation results in a power-grid example that corroborate the technical findings.
\end{abstract}


\maketitle


\section{Introduction}
This paper studies the stability properties of the interconnection 
between a Linear Time Invariant (LTI) dynamical system and an output 
feedback controller, where the control inputs can be maliciously 
compromised by an attacker. 
The theoretical and algorithmic endeavors are motivated by a number of applications within the realm of cyber-physical systems (CPSs) -- that is, physical systems integrated with computational resources by means of a communication infrastructure. While advances in communication and cyber technologies provide enhanced functionality, efficiency, and autonomy of a CPS, the presence of communication channels as well as the tight integration between cyber and physical components unavoidably introduces new security vulnerabilities. In particular,  the modeling adopted in this paper is well suited for a number of applications in power systems~\cite{menta2018stability,MC-ED-AB:19}, transportation systems~\cite{bianchin2020online}, communication networks~\cite{low1999optimization}, and robotics~\cite{zheng2019implicit}, to mention just a few. 

\begin{figure}[pos=t]
\centering 

\includegraphics[width=\columnwidth]{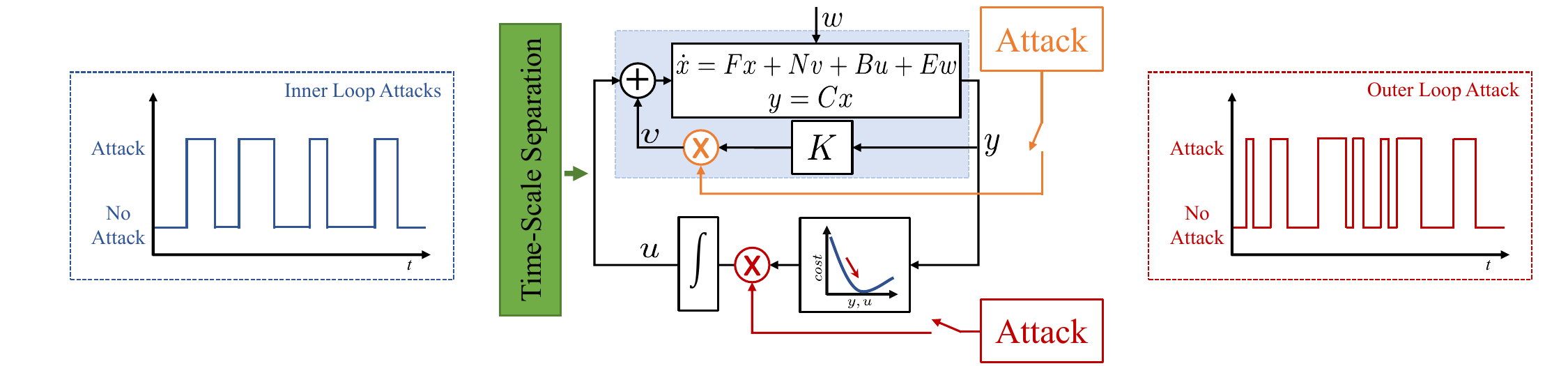}
\caption{An inner control loop is designed to 
stabilize the plant, while an outer control loop regulates the plant 
to an optimizer of \eqref{eq:optpro0}.}
\vspace{-.4cm}
\label{fig:controllerBlocks}
\end{figure}

As illustrated in Figure~\ref{fig:controllerBlocks}, our setting considers two types of feedback controllers: (i) a static output-feedback controller that is used as 
an inner loop to stabilize the LTI system; and, (ii) an output feedback controller based on an appropriate modification of gradient-flow dynamics. For the synthesis of the gradient-flow controller, we start by formulating an unconstrained optimization problem with a composite cost function that captures performance indexes associated with the input and the output of the plant. The cost function is assumed to be smooth and to satisfy the  Polyak-\L ojasiewicz (PL)
inequality condition (where we recall that the PL inequality is a weaker assumption than convexity, and it implies the property of invexity)~\cite{HK-JN-MS:16}. We consider adversarial actions of the form of \textit{switching  multiplicative attacks against the interconnection between the LTI plant and the controller}, whereby an attacker can persistently modify the inputs of the dynamical system with the objective of  destabilizing the equilibrium points. This attack model is rather general and it captures different classes of attacks that can \emph{persistently} modify the sign and/or the  magnitude of the control signals, and also jam the communication 
channels. 

With this setting in place, the problem  addressed in this paper is 
that of finding sufficient conditions that guarantee that the 
equilibrium points of the closed-loop system are asymptotically or 
exponentially stable. To this end, we adopt the framework of hybrid 
systems, and we leverage arguments from singular perturbation 
analysis and input-to-state stability (ISS) for hybrid systems. 
This framework allows us to establish sufficient conditions in terms of the total activation time of the attacks acting on the system, as well as the time-scale separation needed between the plant and the gradient-flow controller to preserve the stability properties of the interconnection.

\textbf{Related Works.} The  design  of feedback-based optimization controllers has recently  received  
significant  attention; see, for example~\cite{Jokic2009controller,brunner2012feedback,hauswirth2020timescale,lawrence2018optimal,MC-ED-AB:19,zheng2019implicit,bianchin2020online,FlorianReview}. In particular, in~\cite{hauswirth2020timescale} sufficient conditions on the time-scale separation between plant and controller were derived to induce asymptotic stability properties. Similarly, the joint stabilization and optimal steady-state regulation of LTI dynamical systems was considered in~\cite{lawrence2018optimal}. LTI dynamical systems with time-varying exogenous inputs were considered in \cite{MC-ED-AB:19}, along with the problem of tracking an optimal solution trajectory of a time-varying problem with a strongly convex cost; prediction-correction-type controllers were utilized to track the trajectory of a time-varying  problem in~\cite{zheng2019implicit}. Recently, \cite{bianchin2020online} established exponential stability results for the interconnection of a switched LTI system and a hybrid feedback controller based on accelerated gradient dynamics with resets. Finally, \cite{FGJ-JIP-GB-ED:21-lcss} and \cite{Poveda:16} studied  extremum-seeking algorithms for model-free optimization under deception attacks. 

In recent years, motivated by the increasing vulnerability of cyber-physical systems operating in adversarial environments, several works have investigated the stability properties of systems under denial-of-service attacks (DoS), whereby an attacker compromises system resources such as
sensors or actuators, as well as infrastructure such as communication channels. We refer to \cite{FP-FD-FB:10y,Senejohnny2018TCN,Cardenas2008} (see also 
references therein) for a comprehensive list of related works, while we present a list of representative references below. The authors in \cite{AC-HI-TH:16} designed a stabilizing controller for communication channels that face malicious random packet losses, while 
\cite{HSF-SM:12} presented a class of stabilizing event-triggered controllers under DoS attacks.
The work  \cite{Kundu2020TCNS} designed scheduling policies that preserve 
the stability of an interconnected system when an attacker jams the control 
inputs that act on the plants.
The works \cite{NF-GB-LC-BS:17} and \cite{SZY-MZ-EF:15} investigated system 
stability in the presence of deception attacks, namely, attacks where the 
integrity of  control packets or measurements is compromised. Robustness of a sliding mode controller against DoS attacks was investigated in~\cite{Chengwei2020DoS}; the effect of DoS attacks in the tracking performance of actuators was  investigated in~\cite{Chenqwei2021Tracking}. Attacks to both a state estimator and actuator of an LTI system were considered in~\cite{Gao2020Attacks}, where the attacks are assumed to satisfy a given sparse observation condition. A self-triggered consensus network in the presence of communication failures caused by DoS attacks was considered in~\cite{Senejohnny2018TCN}, and a similar problem was studied in ~\cite{Wang2020TC} and  \cite{WANG2020AUT}. Event-triggered communication and decentralized control of switched systems under cyber attacks are analyzed in~\cite{QI2020ISA}; in particular, conditions on the dwell-time and the gain for the controller are derived to ensure stability. Other lines of work have studied the robustness properties with respect to attacks of different (discrete-time) optimization algorithms; e.g.,  \cite{SS-BG:18,LS-NHV:20,YC-LS-JX:17}. In these works, Byzantine attacks in distributed algorithms are modeled by a group of malicious nodes that modify the data transmitted to their neighbors. A similar model was considered in~\cite{BT-CAU-HW-MA:20} for sub-gradient
methods. Furthermore, a model-free moving target defense framework for the detection and mitigation of sensor and/or actuator attacks with discrete-time dynamics was considered in~\cite{ZHAI2021104826}.

\textbf{Contributions.} 
The contribution of this work is threefold. 
i) We formulate a new class of switching-mode multiplicative attacks against the feedback signals, whereby the attacker transforms the control inputs produced by the controller according to a linear map with the objective of destabilizing the closed-loop system. This class of attacks is novel in the literature and it includes, as a special case, DoS and deception attacks. ii) We present the first stability analysis of the interconnection between an LTI dynamical system and static or dynamic controllers operating under switching multiplicative attacks. Our framework leverages analytical tools from set-valued hybrid dynamical systems theory, optimization, and feedback control theory. We characterize sufficient conditions on the frequency of the attacks to guarantee exponential stability of the closed-loop system under a suitable time-scale separation between the dynamics of the plant and the dynamics of the controller. iii) For systems with time-varying exogenous inputs acting on the plant dynamics, we establish sufficient conditions to guarantee exponential input-to-state stability with respect to the essential supremum of the time derivative of the exogenous signal, similar to the notion of derivative ISS studied in \blue{ \cite{Angeli2003_diss}}. To derive this result, we leverage Lyapunov theory for switched systems with inputs and unstable modes, where the switching signals are generated by a hybrid dynamical system. 

\textbf{Organization.} The rest of this paper is organized as follows. In Section \ref{eq:secPreliminaries} we present the notation used throughout the paper and some preliminaries on hybrid dynamical systems. In Section \ref{section_problem_formulation}, we formalize the model of the system and the problem under study. Section \ref{section_main_results} presents the main results, followed by the analysis presented in Section \ref{section_analysis}. Section \ref{section_numerical_examples} presents some numerical experiments, and Section \ref{section_conclusions} concludes the paper by summarizing our main findings.

\section{Preliminaries and Notation}
\label{eq:secPreliminaries}
We begin by introducing the notation that will be used throughout the 
paper. Given a compact set $\mathcal{A}\subset\R^n$, 
$n \in \mathbb{Z}_{\geq0}$, and a vector $z\in\R^n$, we let 
$|z|_{\mathcal{A}}:=\min_{s\in\mathcal{A}}\|z-s\|_2$ denote the 
minimum distance between $z$ and $\mathcal{A}$. 
\blue{We denote by $\overline{\mc A}$ the closure of set  $\mc A$.}
We denote $\mathbb{B}:= \{x\in \R^n ~:~\norm{x}\leq 1\}$. For a bounded
function $w: \R_{\geq 0} \rightarrow \R^n$, we 
define the sup norm as $\norm{w}_t := \sup_{s \in [0,t]}\norm{w(s)}$.
For a matrix $A \in \R^{n \times d}$, with
$d \in \mathbb{Z}_{\geq0}$, we use $\norm{A}$ to denote the \blue{largest singular value} of 
$A$. When $n=d$, and the matrix is symmetric, we use $\bar{\lambda}(A)$ and $\underline{\lambda}(A)$ to denote the largest and smallest eigenvalue of $A$, respectively. Given a set $O$, we use $I_O(x)$ to denote the indicator function for the set 
$O$, namely $I_O(x)=1$ if $x \in O$ and
$I_O(x)=0$ if $x\notin O$. Given vectors 
$p_{1}\in \real^d,p_2\in\R^n$, we denote by  $(p_1,p_2) \in \real^{d+n}$ their concatenation. We denote the domain of a function $p$ as $\text{dom}(p)$. A function $\alpha: \R_{\geq 0}\rightarrow \R_{\geq 0}$ is of class-$\mathcal{K}$ if it is continuous, strictly increasing, and $\alpha(0)=0$. If, in addition, $\lim_{r \rightarrow \infty}\alpha (r) = \infty$, it is said to be of class-$\mathcal{K}_\infty$. A function $\sigma:\R_{\geq 0}\rightarrow \R_{\geq 0}$ is of class-$\mathcal{L}$ if it is continuous, decreasing and $\lim_{r \rightarrow \infty}\sigma(r) = 0$. A function $\beta: \R_{\geq 0} \times \R_{\geq 0} \rightarrow \R_{\geq 0}$ is said to be of class-$\mathcal{K}\mathcal{L}$ if it is class-$\mathcal{K}$ in its first argument and class-$\mathcal{L}$ in its second argument.

In this paper, we will consider mathematical models corresponding to set-valued dynamical systems that combine continuous-time 
dynamics (with inputs) and discrete-time dynamics. These systems are referred to as Hybrid Dynamical 
Systems (HDS) \cite{Goebel:12}, and they are characterized by the following hybrid inclusions:
\begin{align}\label{HDS}
p\in \blue{\mc C},~~\dot{p}\in F(p,q),
~~~~~~~~~~p\in \blue{\mc D},~~p^+\in G(p),
\end{align}
where $p\in\mathbb{R}^n$ is the state of the system, 
$q \in \R^m$ is the input, 
$F:\mathbb{R}^n\times \R^m\rightrightarrows\mathbb{R}^n$
is a set-valued map (called the \emph{flow map}) that governs the 
continuous-time dynamics of the system when the state belongs to the 
flow set $\blue{\mc C}\subset\mathbb{R}^n$. Similarly, 
$G:\mathbb{R}^n\rightrightarrows\mathbb{R}^n$
is a set-valued map (called the \emph{jump map}) that describes the 
discrete-time dynamics of the system when the state belongs to the 
jump set  $\blue{\mc D}\subset\mathbb{R}^n$.  The tuple 
$\mathcal{H} = \{\blue{\mc C},F,\blue{\mc D},G\}$ completely
system. \blue{A HDS is well-posed if it 
satisfies~\cite[Def. 6.29]{Goebel:12} and it is nominally well-posed 
if it satisfies~\cite[Def. 6.2]{Goebel:12}. We recall that a 
well-posed hybrid system is nominally well-posed.}
Solutions to \eqref{HDS} are defined on \emph{hybrid time domains}, 
namely, they are indexed by a parameter 
$t\in\mathbb{R}_{\geq0}$ that increases continuously during flow 
(i.e. while the continuous dynamics are active), and by a parameter 
$j\in\mathbb{Z}_{\geq0}$ that increases by one unit during the jumps 
(i.e., when the discrete dynamics are active).  
A set $E\subset \R_{\geq0}\times\mathbb{Z}_{\geq0}$ is called a 
\textsl{compact} hybrid time domain if  
$E=\cup_{j=0}^{J-1}([t_j,t_{j+1}],j)$ for some  
$0=t_0\leq t_1 \leq \ldots\leq t_{J}$. The set $E$ is a hybrid time 
domain if for all $(T,J)\in E$, the set 
$E\cap([0,T]\times\{0,\ldots,J\})$ is a compact 
hybrid time domain. Using the notion of hybrid time domains we can 
formally introduce the concept of solution to systems of the form 
\eqref{HDS}, which is borrowed from \cite[Ch. 2]{Goebel:12}.
\begin{definition}
A function $p:\text{dom}(p) \to \mathbb{R}^n$ is a hybrid arc if 
$\text{dom}(p)$ is a hybrid time domain and $t\mapsto p(t,j)$ is 
locally absolutely continuous for each $j$ such that the interval 
$I_j:=\{t:(t,j)\in \text{dom}(p)\}$ has nonempty interior. A hybrid 
arc $p$ is a \emph{solution} to \eqref{HDS} with 
$\text{dom}(q)=\text{dom}(p)$ if 
$p(0,0)\in \blue{\overline{\mc C}} \cup \blue{\mc D}$, and the following two 
conditions hold:
\begin{enumerate}
\item For each $j\in\mathbb{Z}_{\geq0}$ such that $I_j$ has nonempty 
interior: $p(t,j)\in \blue{\mc C}$ for all $t\in \text{int}(I_j)$, and 
$\dot{p}(t,j)\in F(p(t,j),q(t,j))$ for almost all $t\in I_j$.
\item For each $(t,j)\in\text{dom}(p)$ such that $(t,j+1)\in \text{dom}(p)$: $p(t,j)\in \blue{\mc D}$, and $p(t,j+1)\in G(p(t,j))$.  
\end{enumerate}
\end{definition}
In this paper, we will always work with HDS that satisfy the following Basic Conditions whenever $q$ is constant:
\begin{enumerate}[a)]
    \item $\blue{\mc C}$ and $\blue{\mc D}$ are closed subset of $\R^n$.
    \item $F(\cdot,q)$ is outer semicontinuous and locally bounded relative to $\blue{\mc C}$, $\blue{\mc C} \subset \text{dom }F$, and $F(p,q)$ is convex for every $p \in \blue{\mc C}$.
    \item $G:\mathbb{R}^n\rightrightarrows\mathbb{R}^n$ is outer semicontinuous and locally bounded relative to $\blue{\mc D}$, and $\blue{\mc D} \subset \text{dom }G$.
\end{enumerate}
On the other hand, when $q$ is not constant we will make it a function $\Phi(\cdot)$ of the state $p$, and in this case the Basic Condition b) should hold for the set-valued map $\hat{F}(\cdot):=F(\cdot,\Phi(\cdot))$. By working with hybrid time domains and nominally well-posed HDS (note that the Basic Conditions are sufficient to obtain well-posedness of the
system~\cite[Thm. 6.8]{Goebel:12}) we can exploit suitable 
graphical convergence notions to establish sequential compactness 
results for the solutions of \eqref{HDS}, e.g., the graphical limit 
of a sub-sequence of solutions is also a solution. Such types of 
properties will be instrumental for the robustness analysis of the 
systems studied in this paper.
\begin{definition}
A hybrid solution $p$ is maximal if there does not exist another solution $\psi$ to $\mathcal{H}$ such that $\text{dom}(p)$ is a proper subset of $\text{dom}(\psi)$, and $p(t,j)=\psi(t,j)$ for all $(t,j)\in\text{dom}(p)$. A maximal hybrid solution is said to be complete if its domain is unbounded. 
\end{definition}
In this paper, we are interested in establishing suitable \emph{convergence and stability} properties for a class of dynamical systems under attack. To this end, the following definitions will be instrumental.
\begin{definition}\label{definition_stability}
Given a compact set $\mathcal{A}\subset \blue{\mc C\cup \mc D}$, 
system 
\eqref{HDS} with $u=0$ is said to render $\mathcal{A}$ uniformly 
globally asymptotically stable (UGAS) if there exists a class 
$\mathcal{K}\mathcal{L}$ function $\beta$ such that every solution of
\eqref{HDS} satisfies  
$|p(t,j)|_{\mathcal{A}}\leq \beta(|p(0,0)|_{\mathcal{A}},t+j)$ for 
all $(t,j)\in\text{dom}(p)$. If, additionally, $\beta(r,s)=c_1re^{-c_2s}$ for 
$c_1,c_2>0$, the system is said to render the set $\mathcal{A}$ 
uniformly globally exponentially stable (UGES).
\end{definition}
Note that the stability notion of Definition 
\ref{definition_stability} implies the standard asymptotic stability properties considered in the literature of continuous-time 
systems $(\blue{\mc D}=\emptyset)$ and discrete-time systems 
$(\blue{\mc C}=\emptyset)$. In particular, recall that a continuous-time LTI 
dynamical system given by $\dot{x} =  Ax$ renders the origin 
$\mathcal{A}=\{0\}$ UGES if and only if $A$ is Hurwitz, i.e., 
Re$\{\lambda(A)\}< 0$, for all eigenvalues $\lambda$ of $A$.

When the input $u$ in \eqref{HDS} is not identically zero, the notion
of \emph{input-to-state stability} can be used to qualitatively 
characterize the effect of the input on the stability properties of 
the system. \blue{We refer to \cite{CAI200947} for a formal 
definition of the sup norm $\norm{q}_{(t,j)}$ on hybrid time 
domains.}
\begin{definition}
For every measurable function $q:\text{dom}(q)\to\mathbb{R}^m$ with $\text{dom}(q)=\text{dom}(p)$, system \eqref{HDS} is said to render the compact set $\mathcal{A}$ input-to-state stable (ISS) with respect to $q$ if there exists a class $\mathcal{K}\mathcal{L}$ function $\beta$, and a class $\mathcal{K}$ function $\gamma$, such that every solution of the system satisfies the bound
\begin{equation}
|p(t,j)|_{\mathcal{A}}\leq \beta(|p(0,0)|_{\mathcal{A}},t+j)+\gamma\left(\norm{q}_{(t,j)}\right),
\end{equation}
for all $(t,j)\in\text{dom}(p)$. When $\beta$ has an exponential form, we say that system \eqref{HDS} renders the set $\mathcal{A}$ exponentially input-to-state stable (E-ISS).
\end{definition}
The properties of UGAS, UGES, and ISS for hybrid systems can be readily established via suitable Lyapunov-based conditions; see \cite[Ch.3]{Goebel:12} for sufficient conditions that certify UGAS, \cite[Thm. 1]{TE-FO-ZA-2013} for UGES, and \cite[Thm. 3.1]{CAI200947} for ISS. 
%
\section{Problem Formulation}
\label{section_problem_formulation}
\blue{
In this section, we formalize the problem studied in this work.
We begin by describing the system and controller model in the absence
of attacks. Subsequently, we present our attack model, and
we formalize the joint control and optimization problem that is 
the focus of this work.}

\subsection{Nominal Model of the Plant and the Controller}
\vspace{0.1cm}
We consider plants modeled as LTI dynamical systems of the form:
\begin{equation}\label{eq:plantModel}
\dot{x} = F x + N v + B u + E w,  \quad y = C x,
\end{equation}
where $x \in \R^{n}$ is the state, $u \in \R^{m}$ and 
$v \in \R^{m_v}$ are control inputs, $w \in \R^{q}$ is a time-varying
\emph{unknown} exogenous signal, $y \in \R^p$ is the measurable 
output, and $F,N, B, E, C$ are matrices of suitable dimensions. 
\blue{We consider general cases where the matrix $F$ may not be 
Hurwitz-stable.}
We make the following regularity assumption on the exogenous
input $w$. 

\begin{assumption}\label{as:cts_w}
The exogenous input $w$ is generated by an exosystem of the form $\dot{w}=\Pi(w)$, with $w(0)\in \mathcal{W}$, where
$\Pi(\cdot)$ is Lipschitz continuous, and $\mathcal{W} \subset \bar{w}\mathbb{B}$ is a compact set. Moreover, this exosystem renders the set $\mathcal{W}$ forward invariant. 
\end{assumption}
Throughout this paper, functions $w$ satisfying the 
conditions of Assumption \ref{as:cts_w} are said to be 
\emph{admissible.}

\vspace{0.1cm}
As graphically illustrated in the center panel of Figure 
\ref{fig:controllerBlocks}, we consider two control loops applied to the plant \eqref{eq:plantModel}. First, an inner control loop, of the form of a static output feedback law, is applied via the control input $v$:
\begin{align}\label{eq:stateFeedback}
v = Ky, \quad K \in \R^{m_v \times p}.
\end{align}
Using \eqref{eq:plantModel} and \eqref{eq:stateFeedback}, we denote 
by $A:= F + N K C$ the matrix describing the closed-loop dynamics of 
the plant. \blue{We assume that the inner feedback loop is 
stabilizing, as formalized in the following assumption.}
\begin{assumption}\label{as:stabilityPlant}
For every positive definite matrix $R\in \R^{n \times n}$, there is a unique positive definite matrix $P \in \R^{n \times n}$ such that $A^\top P + P A = -R$. 
\end{assumption}

Additionally, we consider an outer loop of the form of a low-gain 
dynamical controller that acts on the input $u$ as follows (see 
Fig. \ref{fig:controllerBlocks}):
\begin{align}
\label{eq:gradientController}
\dot{u} &=   - \varepsilon \left(\nabla f_u(u) 
+ G^\top \nabla f_y(y)\right),
\quad G:= -C A^{-1}B,
\end{align}
where $f_u: \R^m \rightarrow \R$ and $f_y:\R^p \rightarrow \R$ are 
cost functions (to be characterized below), and $\varepsilon > 0$ is a tunable gain.

By combining the plant model \eqref{eq:plantModel} with the 
controllers \eqref{eq:stateFeedback}-\eqref{eq:gradientController}, 
we obtain the following nominal closed-loop system with states 
$(x,u) \in \R^n \times \R^m$:
\begin{align}\label{eq:interconnection}
\dot{x} &= (F+NKC) x + B u + E w,  \quad y = C x, 
\nonumber\\
\dot{u} &=   -\varepsilon\Big(\nabla f_u(u) + G^\top \nabla f_y(y)\Big).
\end{align}

\blue{As shown in \cite{bianchin2020online,hauswirth2020timescale}, under Assumption 
\ref{as:stabilityPlant} and suitable conditions on the gradients $(\nabla f_u,\nabla f_y)$, there exists a sufficiently small gain $\varepsilon$ such that the controller 
\eqref{eq:gradientController} regulates system 
\eqref{eq:plantModel} towards the solutions of:}
\begin{align}\label{eq:optpro0}
\min_{u,\blue{x},y} ~f_u(u) + f_y(y),~~~~~~~\text{s.t.}~~~~0 = A x + Bu + Ew, ~~~~~~~y = Cx.
\end{align}
where $\map{f_u}{\real^m}{\real}$ is interpreted as a cost function 
associated with the steady-state control input, and 
$\map{f_y}{\real^p}{\real}$ is interpreted as a
cost associated with the steady-state system output.
\blue{
The optimization problem \eqref{eq:optpro0} formalizes an optimal 
equilibrium-selection problem, where the objective is to select 
an input-state-output triplet $(u,x,y)$ that minimizes  the cost 
function $f_u(u) + f_y(y)$ at steady state. Moreover, since \eqref{eq:optpro0} is parametrized by the 
(deterministic) time-varying signal $w$, solutions of \eqref{eq:optpro0}, denoted by $(u^*, x^*, y^*)$ will also be time-varying.
}
\begin{remark}
Steady state optimization problems of the form \eqref{eq:optpro0} have received increased attention in the literature during the last years, partially motivated by applications in power 
systems \cite{MC-ED-AB:19,hauswirth2018time,lawrence2018optimal,hauswirth2020timescale} and transportation networks 
\cite{GB-JC-JP-ED:21-tcns}. Note that when $A$ is Hurwitz the steady-state constraint can be solved for $x$ and substituted in the cost function, leading to an unconstrained optimization problem defined over $u$.
\end{remark}
Let $G:= -C A^{-1}B$ and $H:=-CA^{-1}E$, where the inverse of $A$ is well-defined thanks to Assumption \ref{as:stabilityPlant}. Problem \eqref{eq:optpro0} can be written as: 
%
\begin{align}\label{eq:optpro00}
    \min_u f(u,w) := f_u(u) + f_y(Gu + Hw),
\end{align}
where $w$ acts as an exogenous signal that parametrizes the optimal 
solution. 
Note that a solution of~\eqref{eq:optpro0} is also a 
solution of \eqref{eq:optpro00}, but the reverse implication may
not be true, unless the pair $ (C,A) $ is observable.
Nonetheless, observability is not needed for the subsequent analysis,
since we focus on solving problem \eqref{eq:optpro00}. \blue{
In the sequel, we use the following expressions for the partial 
gradients of the costs:
\begin{align*}
    \nabla_u f(u,w) = \nabla_u f_u(u) + G^\top\nabla_y f_y(Gu + Hw), \qquad \nabla_{w} f(u,w) = H^\top \nabla_y f_y(Gu + Hw).
\end{align*}}
We make the following standard regularity assumptions:
\begin{assumption}\label{as:lipschitzcost}
The functions $u \mapsto f_u(u)$ and $y \mapsto f_y(y)$ are 
continuously differentiable and globally Lipschitz, namely, there exist
$\ell_u,~\ell_y >0$ such that 
$\norm{\nabla f_u(u) - \nabla f_u(u')} \leq \ell_u \norm{u - u'}$ and
$\norm{\nabla f_y(y) - \nabla f_y(y')} \leq \ell_y \norm{y - y'}$
hold for every $u, u' \in \R^m$ and $y, y' \in \R^p$.
\end{assumption}
\begin{assumption}\label{as:PLineq}
For each $w$, \blue{the function $u \mapsto f(u,w)$ has a 
unique (bounded) minimizer, denoted by $u^*_w$}. Furthermore,  the function 
$u \mapsto f(u,w)$ satisfies the PL inequality, uniformly in $w$,
namely, $\exists~ \mu > 0$ such that 
$\frac{1}{2}\norm{\nabla_u f(u,w)}^2 \geq \mu 
(f(u,w) - f(u^*_w,w))$, 
$\forall u \in \R^m, w \in \mathcal{W}$. 
\end{assumption}
\begin{assumption}\label{as:cts_ustar}
For each $w\in\mathcal{W}$, the solution of \eqref{eq:optpro00} satisfies $u^*_w = h(w)$, where $w \mapsto  h(w)$ is a continuously differentiable function. 
\end{assumption}
\blue{
\begin{remark}
Assumption \ref{as:lipschitzcost} guarantees that the cost functions 
are sufficiently smooth, which is a standard assumption in the 
literature of online optimization (see also 
\cite{MC-ED-AB:19,lawrence2018optimal}). Indeed, 
by Assumption \ref{as:lipschitzcost}, the mapping
$u \mapsto f(u,w)$ has a globally Lipschitz-continuous gradient, 
uniformly in $w$, with Lipschitz constant 
$\ell = \ell_v + \ell_y ||G||^2$. On the other hand, Assumption \ref{as:PLineq} implies that $f$ is invex, and it is one of the weakest assumptions in the optimization literature that ensures exponential convergence to a set of global optimizers \cite{HK-JN-MS:16}. We also  note that the PL inequality implies the quadratic growth condition $f(u,w) - f(u_w^*,w) \geq \frac{\mu}{2}\norm{ u - u_w^*}^2$, 
$\forall u \in \R^m$ for a fixed $w$; since the minimizer is assumed to be unique and bounded (cf. Assumption \ref{as:cts_w}, \ref{as:PLineq} and \ref{as:cts_ustar}), $\mathcal{W}$ is compact, and $h(\cdot)$ is continuous, the quadratic growth condition  implies radial unboundedness of $f$. Finally, note that by Assumption \ref{as:stabilityPlant}, and for fixed $u$ and $w$, there is a unique solution $x^*_w$ to the constraint $0 = Ax +Bu + Ew$. Thus, by Assumptions \ref{as:cts_w}-\ref{as:cts_ustar}, the solution of \eqref{eq:optpro00} exists, it is unique for each time $t$, and it is a continuous function of $w$.
\end{remark}
}
%
%
Under Assumptions \ref{as:lipschitzcost}--\ref{as:PLineq}, the 
following chain of inequalities holds:  
\begin{align*}
    \frac{\mu^2}{2}\norm{u - u^*_w}^2 \leq \mu(f(u,w) - f(u^*_w,w)) \leq \frac{1}{2}\norm{\nabla_u f(u,w)}^2 \leq \frac{\ell^2}{2}\norm{u-u^*_w}^2.
\end{align*}
The following two lemmas will also be instrumental for our results.
\begin{lemma}\label{lm:derivative_w}
Let Assumptions \ref{as:cts_w} and \ref{as:cts_ustar} hold. Then  $\exists~U_{\bar{w}} \in \mathbb{R}_{> 0}$ such that   $ \norm{\nabla_{w} h(w)}\leq U_{\bar{w}}$, for all $w \in \mathcal{W}$.
\end{lemma}
\begin{lemma}
\label{lm:derivative_w2}
Let Assumptions \ref{as:lipschitzcost}-\ref{as:cts_ustar} hold. Then, $\nabla_u f(h(w),w)^\top \nabla_{w}h = 0$.

\end{lemma}
The result of Lemma~\ref{lm:derivative_w} follows from the continuity of $\nabla_{w} h(\cdot)$ (with respect to $w$) and the compactness of the set $\mathcal{W}$. Regarding Lemma~\ref{lm:derivative_w2}, by first-order optimality conditions the gradient of the function $f$ is equal to zero at the set of optimizers; the product is well-defined since $\nabla_u f(u, w)$ and $\nabla_{w} h$ are continuous by Assumptions \ref{as:lipschitzcost} and \ref{as:cts_ustar}, and $\nabla_w h$ is bounded by Lemma \ref{lm:derivative_w2}, hence, the result. 


\subsection{Attack Model}

\vspace{0.1cm}
In this work, we consider attacks that can modify the control signals
generated by the controllers. As also illustrated in 
Figure~\ref{fig:controllerBlocks}, we focus on attacks that are 
intermittent and occur in a persistent fashion, but are not 
necessarily periodic. We account for two types of attacks. First, 
we consider attacks acting on the gradient-flow 
controller~\eqref{eq:gradientController}, which can be modeled as:
\begin{subequations}
\label{eq:attackModelDynamic}
\begin{align}
\text{Nominal System: } & \dot u = - \varepsilon 
\left(\nabla f_u(u) + G^\tsp \nabla f_y(y)\right), \\
\text{System Under Attack: } & \dot u = - \varepsilon  M_{\sigma_u}
\left(\nabla f_u(u) + G^\tsp \nabla f_y(y)\right),
\end{align}
\end{subequations}
where $\map{\sigma_u}{\realpos}{\Sigma_{u,a}}$
is a switching signal taking values in a finite set 
of indices $\Sigma_{u,a}\subset\mathbb{Z}_{>0}$, and 
$M_{\sigma_u}\in \real^{m \times m}$ is a  matrix that describes the 
attack transformation map. \blue{The 
model~\eqref{eq:attackModelDynamic} is general enough to capture 
different types of attacks, including attacks that modify the sign 
and/or the magnitude of (some or all components of) the gradient (see
e.g.~\cite{SS-BG:18,LS-NHV:20,YC-LS-JX:17,BT-CAU-HW-MA:20,Gao2021Aperiodic, Songlin2020Aperiodic,Zhang2019Aperiodic} for similar attack models), jamming attacks, and DoS attacks  (see e.g.~\cite{Senejohnny2018TCN,WANG2020AUT}).}

\vspace{0.1cm}
Second, we consider attacks targeting the static feedback controller \eqref{eq:stateFeedback}, which can be modeled as:
\begin{subequations}
\label{eq:attackModelStatic}
\begin{align} 
\text{Nominal System: }~~ & v = K y , ~~~~~~~~~~ \dot{x} = F x + N v + B u + E w, \\
\text{System Under attack: } ~~& v = L_{\sigma_v} K y , ~~~~ \dot{x} = F x + N v + L^b_{\sigma_v} B u + L^e_{\sigma_v} E w,
\end{align}
\end{subequations}
where $\map{\sigma_v}{\realpos}{\Sigma_{v,a}}$ is a 
switching signal taking values in a finite set of 
indices $\Sigma_{v,a}\subset\mathbb{Z}_{>0}$. For each $\sigma_v  \in \Sigma_{v,a}$, 
$L_{\sigma_v} \in \real^{m_v \times m_v}$ is a matrix that describes the transformation map applied by the attacker to the static feedback loop, and the matrices 
$L^e_{\sigma_v}, L^b_{\sigma_v} \in \real^{n \times n}$ modify the 
inputs $B u$ and $E w$, respectively.

\blue{
For our analysis, we next unify the models of the 
\emph{Nominal System} and of the \emph{System Under Attack} by 
combining the two modes into a single switching-system model.
Tot his aim, we describe the \emph{Nominal System} in 
\eqref{eq:attackModelStatic} by introducing the index set 
$\Sigma_{v,s}=\{s\}$ and the identity transformation map
$L_s=I \in \real^{m_v \times m_v}$.
Similarly, we describe the \emph{Nominal System} in 
\eqref{eq:attackModelDynamic} by introducing the index set  
$\Sigma_{u,s}=\{s\}$ and the identity transformation map
$M_s=I \in \real^{m \times m}$. }
Using these definitions, we can rewrite systems \eqref{eq:attackModelDynamic}-\eqref{eq:attackModelStatic} as switching systems of the form
\begin{align}\label{eq:attackModel}
\dot u &= - \varepsilon  M_{\sigma_u} \left(\nabla f_u(u) + G^\tsp \nabla f_y(y)\right), \\ 
v &= L_{\sigma_v} K y,  \hspace{1.0cm} B_{\sigma_v} :=   L^b_{\sigma_v} B, \hspace{1.0cm}  E_{\sigma_v} :=   L^e_{\sigma_v} E,
\end{align}
where now the switching signal $\sigma_v$ takes values in the extended 
set $\Sigma_v:=\Sigma_{v,s} \cup \Sigma_{v,a}$, and the switching signal $\sigma_u$ takes values in the extended set 
$\Sigma_u:= \Sigma_{u,s} \cup \Sigma_{u,a}$.

\begin{assumption}\label{as:attackModel}
The sets $\Sigma_v$ and $\Sigma_u$ are finite.
Additionally,  any attack is destabilizing, namely:
\begin{enumerate}
\item[(i)] For all $\sigma_v \in \Sigma_{v,a}$, the matrix $F+N L_{\sigma_v} K C$ is invertible and admits at least one eigenvalue with positive real part.
\item[(ii)] For all $\sigma_u \in \Sigma_{u,a}$, the matrix $M_{\sigma_u}$ has at least 
one eigenvalue with negative real part. Moreover, we define $ \bar{M}:=  \max_{\sigma_u \in \Sigma_{u,a}} \norm{M_{\sigma_u}}$.

\end{enumerate}
\end{assumption}
It follows from Assumption \ref{as:attackModel} that modes 
described by the sets $\Sigma_{v,a}$ and $\Sigma_{u,a}$ can make the closed-loop system unstable \eqref{eq:interconnection}. Moreover, 
since the sets $\Sigma_{v,a}$ and $\Sigma_{u,a}$ are countable and  
finite, with a slight abuse of notation, in the remainder we will 
enumerate these unstable modes as follows:
$\Sigma_{v,a} := \{a_1, \dots ,a_{\vert \Sigma_{v,a} \vert}\}$ 
and $\Sigma_{u,a} := \{a_1, \dots ,a_{\vert \Sigma_{u,a} \vert}\}$, 
where the distinction between the two sets will be made clear by the 
context.

\begin{lemma}\label{lm:symR}
Let Assumption \ref{as:attackModel}-(i) hold. Let $A_{a_i}: = F + NL_{a_i}KC$ for  
$a_i \in \Sigma_{v,a}$, and let  $\hat{R}_{a_i} := (A_{a_i}^\top  +  A_{a_i})$. Then, $\exists ~\lambda_i >0,$ and $x_i \in \mathbb{C}^n$, $x_i \neq 0$ such that $\hat{R}_{a_i}x_i = \lambda_i x_i$\footnote{The existence of eigenvectors and eigenvalues in this case can be proved as in \cite{existenceEigenvalues}.}.
\end{lemma}

\emph{Proof: } Suppose by contradiction that $\hat{R}_{a_i} \preceq 0$, and let $x_a \in \mathbb{C}^n$ be a non-zero vector such that $A_{a_i} x_a =  \lambda_a x_a$, then\footnote{We denote the conjugate of a matrix $x \in \mathbb{C}^{m \times n}$ as $x^H$.}: 
\begin{align*}
x_a^H \hat{R}_{a_i} x_a &
    = x_a^H (A_{a_i}^\top + A_{a_i})x_a 
    = (x_a^H \lambda_a^H x_a + x_a^H  \lambda_a x_a)  
    = (\lambda_a^H + \lambda_a)x_a^H x_a
    = 2 \mathrm{Re}(\lambda_a)x_a^H x_a,
\end{align*}
and since $\hat{R}_{a_i} \preceq 0$ and $x_a\neq0$, we must have that  $\mathrm{Re}(\lambda_a)\leq 0$. Since $\lambda_a$ is arbitrary, this contradicts 
Assumption \ref{as:attackModel}-(i).  \hfill $\blacksquare$\\

\noindent
For any $a_i \in \Sigma_{v,a}$, the matrix $\hat{R}_{a_i}$ is real and symmetric; hence, Rayleigh theorem holds \cite[Thm. 4.2.2]{horn_johnson_2012} and  $\underline{\lambda}(\hat{R}_{a_i}) \norm{x}^2 \leq x^\top \hat{R}_{a_i} x \leq \bar{\lambda}(\hat{R}_{a_i})\norm{x}^2$. 
Throughout the rest of this paper, we denote by $\bar{\lambda}(\hat{R}_a) := \max_{a_i \in \Sigma_{v,a}} \bar{\lambda}(\hat{R}_{a_i})$ the largest eigenvalue of all matrices $\hat{R}_{a_i}$.

\subsection{Admissible Persistent Attacks}
\vspace{0.1cm}

The models \eqref{eq:attackModelDynamic} and \eqref{eq:attackModelStatic} describe a family of persistent multiplicative attacks that intermittently modify the control signal 
generated by the controllers \eqref{eq:stateFeedback} and \eqref{eq:gradientController} via the switching maps $L_{\sigma_v}$  and $M_{\sigma_u}$. As described in Assumption 
\ref{as:attackModel}, such attacks can induce instability in the closed-loop system \eqref{eq:interconnection} by persistently altering the direction of the vector field in the dynamic 
controller \eqref{eq:gradientController}, or the gains of the feedback law \eqref{eq:stateFeedback}. 
We note that, because the attack is multiplicative rather than additive, in general, standard approaches used to guarantee closed-loop stability under non-persistent adversarial attacks (such as 
dynamic filters \cite{NF-GB-LC-BS:17}, identification mechanisms \cite{SZY-MZ-EF:15}, dynamic compensators \cite{HSF-SM:12}, etc) may not guarantee that the attack is mitigated at all times.
For this reason, in this work, we investigate cases where such destabilizing attacks can be tolerated during bounded fractions of time. Specifically, for every time-interval $0 \leq s < t$, we define the total activation time of the attacks \eqref{eq:attackModelDynamic}, and \eqref{eq:attackModelStatic}, respectively, as follows:
\begin{align*}
T_u(s,t):= \int_s^t  I_{\Sigma_{u,a}}(\sigma_u(\tau)) d\tau,~~~~~~~\text{and}~~~~~~~~~T_v(s,t):= \int_s^t I_{\Sigma_{v,a}}(\sigma_v(\tau)) d\tau,
\end{align*}
where we recall that $I_\Sigma(\sigma)$ denotes the indicator 
function for the set $\Sigma$ (see Section 
\ref{eq:secPreliminaries}).
Notice that, if $\sigma_v = \sigma_u = s$ at all times, then 
$T_v(s,t)=T_u(s,t)=0$ for any $0\leq s<t$ since no attack is active 
in the closed-loop. 
The class of admissible attacks is formalized next.

\begin{definition}{\bf \textit{(Class of Admissible Attacks)}} 
\label{def:undetectableAttack}
Given $\kappa_{v,1}, \kappa_{u,1}, \kappa_{v,2}, \kappa_{u,2} >0$ and $N_{0,v},N_{0,u},T_{0,v},T_{0,u} \in \integ_{>0}$, the attacks \eqref{eq:attackModelDynamic} and \eqref{eq:attackModelStatic} are said to be admissible for system \eqref{eq:interconnection} if the number of discontinuities of $\sigma_v$ and $\sigma_u$ in the time-interval $[s,t)$, denoted by $N_v(s,t)$ and $N_u(s,t)$, respectively, satisfy
 \begin{align}\label{eq:avgdwelltimeCondition}
 (a)\qquad
    N_v(s,t) \leq \kappa_{v,1}(t-s) + N_{0,v}, \quad
    N_u(s,t) \leq \kappa_{u,1}(t-s) + N_{0,u},
\end{align}
and the activation times $T_v(s,t)$ and $T_u(s,t)$ satisfy:
\begin{align}
\label{eq:activationTimeCondition}
 (b)\qquad
T_v(s,t) \leq \kappa_{v,2} (t-s) + T_{0,v}, \quad 
T_u(s,t) \leq \kappa_{u,2} (t-s) + T_{0,u}.
\end{align}
\end{definition}
Condition \eqref{eq:avgdwelltimeCondition} imposes an average dwell-time constraint on the switches of $\sigma_u$ and $\sigma_v$, and thus guarantees that any admissible attack cannot 
generate Zeno behavior. On the other hand, condition \eqref{eq:activationTimeCondition} says that the total activation
time of any admissible persistent attack in a given time interval must grow at most linearly with the length of the time interval. We note that larger values of $\kappa_{v,2}$ and $\kappa_{u,2}$ describe more intrusive attacks, namely, attacks the can remain active for longer periods of time. Finally, we observe that similar attack models have been considered in the literature (see e.g. \cite{Poveda:16,Senejohnny2018TCN, WANG2020AUT,FGJ-JIP-GB-ED:21-lcss}). Indeed, conditions \eqref{eq:avgdwelltimeCondition}-\eqref{eq:activationTimeCondition} are rather general as they capture periodic or aperiodic attacks, and they also allow a finite number of consecutive switches in every interval of time of sufficiently large length. 

\begin{remark}
\blue{In cases where the closed-loop system is monitored by an 
attack-detection mechanism, the quantities 
$\{T_{0,v},\kappa_{v,2},T_{0,u},\kappa_{u,2} \}$ can be interpreted 
as attack design parameters that can be tuned by the attacker to
remain undetected from the monitor.}
Alternatively, \blue{in cases where the closed-loop system includes 
an 
attack-mitigation mechanism that can reject the effects of attacks 
for limited periods of time, the quantities 
$\{T_{0,v},\kappa_{v,2},T_{0,u},\kappa_{u,2} \}$ can be 
interpreted as parameters that describe the effectiveness of 
the mitigation mechanism.} 
Both engineering interpretations will be consistent with
our theoretical framework.
\end{remark}

\subsection{Control Objectives: Approximate Optimal Tracking Under Attacks}

\vspace{0.1cm}
Based on the model introduced in the previous sections, we can now formalize the problem under study in this paper.  In particular, for the family of closed-loop systems under attacks described by system \eqref{eq:interconnection}, and satisfying Assumptions \ref{as:stabilityPlant}-\ref{as:attackModel}, we are interested in characterizing the \blue{admissible persistent attacks}, and the time-scale separation between the plant and the controller needed to preserve the stability properties of the nominal closed-loop system. Moreover, when $w$ is constant, we are interested in guaranteeing exponential convergence of the trajectories of the closed-loop system to the set of solutions of the online optimization problem described by \eqref{eq:optpro0}. 
Notice that, when $w$ is time-varying, convergence to the solution 
of the optimization problem can only be achieved up to an error 
that depends on $\dot{w}$. In this case, we seek to quantify this error using the notion of input-to-state stability.

Based on this, the following problem formalizes the stated objectives
\begin{problem}\label{prob:controlObjective}
Determine, when possible, a set of parameters $(\varepsilon, \kappa_{v,1},\kappa_{v,2})$, or $(\varepsilon,\kappa_{u,1} \kappa_{u,2})$ such that the optimal point $z^* = (x^*,u_w^*)$ of system \eqref{eq:interconnection} under attack is exponentially input-to-state stable (E-ISS) with respect to the time-variation (i.e. the derivative) of the exogenous signal  $t\mapsto w(t)$.
\end{problem}
%
%
\begin{remark}
For simplicity and space limitations, we will focus on attacks that are exclusive either to the plant or to the controller, i.e. we do not address simultaneous attacks to the static and dynamic feedback controller. This allows us to simplify our notation since in this case, if $\sigma_v = a_i$, then $\sigma_u = s$, and, if  $\sigma_u = a_i$, then $\sigma_v = s$ as defined before. Note, however, that considering simultaneous attacks in the plant and the controller is a natural extension of this work that can be studied following similar steps as in this paper.
\end{remark}

Before presenting our main results, we list three representative examples of applications in cyber-physical systems for which the modeling and control framework considered in this paper is particularly well-suited. 
\begin{example} The dynamical model~\eqref{eq:plantModel} and the interconnection~\eqref{eq:interconnection} fit within the context of real-time frequency control and economic optimization of power transmission systems. In particular, the matrix $K$ can be designed based on the so-called automatic generation control (AGC) as well as pertinent droop controllers implemented in the generation units~\cite{simpson2020area}; on the other hand, the gradient-flow controller in~\eqref{eq:interconnection} can be utilized to produce setpoints for the generators (both conventional fossil-fuel and renewable-based units) to solve an economic dispatch (ED) or optimal power flow (OPF) problem in real time~\cite{chakraborty2020dynamics}; see also the example provided in~\cite{MC-ED-AB:19}. The models of attacks considered in this paper capture in this case malicious attacks to the AGC signals, droop control loops, and ED/OPF commands.  
\end{example}
\begin{example}
The model ~\eqref{eq:plantModel} and~\eqref{eq:interconnection} is also suitable for applications in the context of intelligent transportation systems, where the properties of safety and resilience are critical. For example, for a ramp metering problem in a highway system, our inner control loop can model local controllers such as the ALINEA~\cite{MP-AK:02}. On the other hand, gradient-flow controllers have been used to drive the equilibrium point of the traffic flows on the highway system towards the solution of a network-level problem as explained in, e.g.,~\cite{GB-JC-JP-ED:21-tcns} and~\cite{grandinetti2018distributed}. It then follows that the attack models consider in \eqref{eq:attackModel} can be used to capture adversarial jamming of the feedback gains of the ALINEA controller, or attacks acting on the control directions obtained from the gradient flow algorithm. 
\end{example}
\begin{example}
The model \eqref{eq:attackModel} is also suitable for applications in robotics, such as the control of autonomous vehicles using navigation laws based on potential and anti-potential functions; see \cite{ObstacleConvexPotentials,SourceSeekingSCL,PovedaTAC21Obstacle}. Indeed, in this scenario, it is common to approximate the vehicle dynamics as a linear system \eqref{eq:plantModel} for which an internal feedback loop is designed to stabilize an external reference $u$; see \cite{BookMurray}. In general, the internal controller is designed to guarantee a unitary DC gain, and to operate in a faster time scale (limited by the actuator dynamics of the vehicle), such that the ``steady state'' of the vehicle satisfies $x^*=u$, where $x$ can model positions or velocities. The input $u$ is then regulated via a gradient flow of the form \eqref{eq:attackModel} aiming to converge to the minimizer (e.g., the target point) of an artificial potential field. In this case, to prevent the vehicle from converging to the target, an external attacker can persistently alter the sign of the gradient dynamics, $i.e., M_{\sigma_u}\in \{I,-I\}$, effectively moving the vehicle \emph{away} from the target. Note that this type of attacks will be particularly detrimental for navigation laws that implement anti-potential fields to avoid obstacles, given that in this case, the attacks will effectively push the vehicles \emph{towards} the obstacles.
\end{example}

\section{Main Results}
\label{section_main_results}

To address Problem \ref{prob:controlObjective}, we will use the framework of HDS \eqref{HDS}, which allows us to model the switching systems \eqref{eq:attackModelDynamic}-\eqref{eq:attackModelStatic} as time-invariant (hybrid) systems, and also to consider suitable (set-valued) dynamic models of the admissible persistent attacks. For simplicity of exposition, we drop the subscripts from \eqref{eq:avgdwelltimeCondition} - \eqref{eq:activationTimeCondition}, 
 and we consider a hybrid automaton with state $\varrho: = (\tau_1,\tau_2,\sigma) \in  \R_{\geq 0} \times \R_{\geq 0} \times \Sigma$, and the following hybrid dynamics:
\begin{subequations}\label{eq:automaton}
	\begin{align}
		\varrho \in \mathcal{C}_\varrho &:=  [0,N_0] \times [0,T_0] \times \Sigma,\\
		\begin{pmatrix} \dot{\tau}_1 \\ \dot{\tau}_2 \\\dot{\sigma}  \end{pmatrix} &\in F(\varrho) := \begin{pmatrix}
		    [0,\kappa_1]\\
			[0,\kappa_2] - I_{\Sigma_a} (\sigma)\\
			0
		\end{pmatrix},\\
		\varrho \in \mathcal{D}_\varrho &:=  [1,N_0] \times [0,T_0] \times \Sigma,\\
		\begin{pmatrix} \tau_1^+\\ \tau_2^+ \\ \sigma^+  \end{pmatrix} &\in G(\varrho) := \begin{pmatrix}
			\tau_1 - 1\\
			\tau_2\\
			\Sigma \backslash \sigma
		\end{pmatrix},
	\end{align}
\end{subequations}
where $T_0\geq 0$, $N_0 \in \mathbb{Z}_{\geq 1}$, $\kappa_1>0$, $\kappa_2 \in (0,1)$. The following Lemma, originally established in \cite{Poveda:16}, will be instrumental for our results.

\begin{lemma}\label{monitoringsystem}
Every complete solution $\varrho$ of \eqref{eq:automaton} has a hybrid time-domain dom$(\varrho)$  that satisfies condition \eqref{eq:avgdwelltimeCondition} and \eqref{eq:activationTimeCondition} with $N_0 = N_{0,v} = N_{0,u}$, and $T_0 = T_{0,v} = T_{0,u}$. Moreover, every signal $\sigma$ modeling an admissible attack can be generated by the dynamical system \eqref{eq:automaton} with a suitable initial condition.
\end{lemma}
Lemma \ref{monitoringsystem}  allows us to study the closed-loop system under persistent switching attacks by studying the interconnection between the dynamics \eqref{eq:attackModelDynamic}, \eqref{eq:attackModelStatic} and \eqref{eq:automaton}. Similar models have been studied in the literature for the analysis of \emph{static} model-based and model-free optimization problems; see \cite{FGJ-JIP-GB-ED:21-lcss,Poveda:16}. However, unlike these works, Problem \ref{prob:controlObjective} describes an \emph{online optimization problem} with a \emph{dynamic} plant in the loop, which has not been studied before in the literature.

In the ensuing section, we derive sufficient conditions for the stability of the interconnected system when the gradient-flow controller is affected by attacks as in \eqref{eq:attackModelDynamic}. After this, we establish similar results for the case when the static feedback controller is under attack, as in \eqref{eq:attackModelStatic}.  
\subsection{Attacks Against Dynamic Feedback}
\vspace{0.1cm}
In this section, we analyze the stability properties of the closed-loop system with attacks acting on the dynamic gradient-based controller. To model this scenario, we consider the following HDS with state $\vartheta:=(x,u,\tau_{u,1},\tau_{u,2},\sigma_u,w)$, and dynamics:
\blue{
\begin{subequations}\label{eq:intercon_u}
	\begin{align}
		\vartheta \in \mathcal{C} &:=  \R^n \times \R^m \times [0,N_{0,u}] \times [0,T_{0,u}] \times \Sigma_u \times \mathcal{W},\\
		\begin{pmatrix} \dot{x}\\\dot{u} \\\dot{\tau}_{u,1}\\
		\dot{\tau}_{u,2} \\ 
		\dot{\sigma}_u \\
		\dot{w}
		\end{pmatrix} &\in F_{\sigma_u}(\vartheta) := \begin{pmatrix}
			A x + B u + Ew\\
			- \varepsilon M_{\sigma_u}\left(\nabla f_u(u) + G^\top \nabla f_y(y)\right) \\
			 [0,\varepsilon\kappa_{u,1}]\\
			[0,\varepsilon\kappa_{u,2}] - \varepsilon I_{\Sigma_{u,a}} (\sigma_u)\\
			0\\
			\Pi(w)
		\end{pmatrix}, \quad y = Cx,\label{flowsmodel1}\\
		\vartheta \in \mathcal{D} &:=  \R^n \times \R^m \times [1,N_{0,u}]\times [0,T_{0,u}] \times \Sigma_u \times \mathcal{W},\\
		\begin{pmatrix} x^+\\ 
		u^+\\
		\tau_{u,1}\\
		\tau_{u,2}^+ \\ 
	    {\sigma_u}^+ \\
	    w^+
	    \end{pmatrix} &\in G_{\sigma_u}(\vartheta) := \begin{pmatrix}
			x\\
			u\\
			\tau_{u,1}-1\\
			\tau_{u,2}\\
			\Sigma_u \backslash \{\sigma_u\}\\
			w
		\end{pmatrix}.
	\end{align}
\end{subequations}
} This HDS models the interconnection between the signal generator \eqref{eq:automaton}, the closed-loop system \eqref{eq:interconnection} with the dynamic controller under attack \eqref{eq:attackModelDynamic}, and the exosystem that generates $w$, which is treated as an auxiliary state. Note that this is similar to the setting considered in \cite[Prop. 5.1]{Angeli2003_diss} in the context of derivative ISS.
%
\begin{remark}
In \eqref{flowsmodel1}, the dynamics of the states $\tau_{u,1}$ and $\tau_{u,2}$ are multiplied by $\varepsilon$. This is done without loss of generality, and only to simplify the modeling framework such that the attacks operate in the same time scale as the gradient controller. In fact, since in our case $\varepsilon\in(0,1)$, the dynamics of \eqref{flowsmodel1} could also be written in terms of the scaled time ratios $\tilde{\kappa}_{u,1}=\varepsilon\kappa_{u,1}$ and $\tilde{\kappa}_{u,2}=\varepsilon\kappa_{u,2}$, and a  $\varepsilon$-scaled indicator function. 
\end{remark}
The following Lemma, which follows directly by construction, will play an important role in the robustness analysis of system \eqref{eq:intercon_u}.
\begin{lemma}\label{wellposedHDS}
For the HDS \eqref{eq:intercon_u} the following holds: (a) The sets $\mathcal{C}$ and $\mathcal{D}$ are closed; (b) The set-valued mapping $F_{\sigma_u}(\cdot) $ is outer-semicontinuous, locally bounded, and convex-valued in $\mathcal{C}$; (c) The set-valued mapping $G_{\sigma_u}$ is outer-semicontinuous and locally bounded in $\mathcal{D}$.  
\end{lemma}
\vspace{0.1cm}
For system \eqref{eq:intercon_u}, we are interested in characterizing sufficient conditions on the gain $\varepsilon$ of the controller, the dwell-time parameter $\kappa_{u,1}$, and the time-ratio parameter $\kappa_{u,2}$, to guarantee the solution of Problem 1. To do this, we first neglect the plant dynamics by assuming they are instantaneous. Namely, we disregard the dynamics $\dot{x}$, and we substitute $x$ by the steady-state mapping $x(u,w)=A^{-1}Bu+Ew$, where $w$ is assumed to be a constant, i.e., $\dot{w}=0$. For this simplified system, we establish the following proposition, which is the first result of this paper. It establishes a novel characterization of an admissible persistent attack.
\begin{proposition}\label{pro:kappa_u2}
Suppose that Assumptions 
\ref{as:cts_w}-\ref{as:attackModel} hold, and that the dwell-time parameter $\kappa_{u,1}$, and the time-ratio parameter $\kappa_{u,2}$ simultaneously satisfy:
\begin{align}\label{timeratiocondition2}
2\mu > \ln(\omega)\kappa_{u,1} + 2\kappa_{u,2}(\mu + \bar{M}\ell), ~~~\text{with}~~~\omega : = \max\left\{\frac{1}{\mu},\frac{\ell^2}{\mu}\right\}.
\end{align}
\blue{
Then, for each compact set $K_x\subset\mathbb{R}^n$, the set 
$$\mc{A}^*_{K_x} =\left\{(x,u,\tau_{u,1},\tau_{u,2},\sigma_u,w):x \in K_x,~u=u_w^*,~\tau_{u,1}\in[0,N_{0,u}],\tau_{u,2}\in [0,T_{0,u}],\sigma_u\in \Sigma_u,w \in \mathcal{W}\right\},$$ 
is UGES for the HDS \eqref{eq:intercon_u} with $\dot{x}=0$, $\varepsilon=1$, $y =  Gu + Hw$, $\dot{w}=0$ and with restricted flow and jump sets 
\begin{align*}
\hspace{-0.5cm}
\mathcal{C}_{K_x} =  (\R^n \cap K_x) \times \R^m \times [0,N_{0,u}] \times [0,T_{0,u}] \times \Sigma_u \times \mathcal{W},~
\mathcal{D}_{K_x} =  (\R^n\cap K_u) \times \R^m \times [1,N_{0,u}] \times [0,T_{0,u}] \times \Sigma_u \times \mathcal{W}. 
\end{align*}
}
\end{proposition}
\vspace{-0.5cm}
The result of Proposition \ref{pro:kappa_u2} establishes an explicit characterization of parameters $(\kappa_{u,1},\kappa_{u,2})$ that guarantees exponential stability of the plant dynamics under persistent attacks. In particular, by \eqref{timeratiocondition2}, we obtain that an admissible attack is dependent on the largest value $\norm{M_{\sigma_u}}$ induced by the attacks in the controller, and the function parameters $\mu$ and $\ell$. Moreover, the parameters $(\kappa_{u,1},\kappa_{u,2})$ are uniform on the initial conditions of the plant and also on the controller. This uniformity property, established via hybrid Lyapunov functions in the next section, will be instrumental for the analysis of the closed-loop system for the case when $\dot{x}\neq 0$.

\begin{remark}
The restriction to a compact set $K_x$ of the $x$-component of the flow and jump set is imposed for the purpose of regularity --namely, to guarantee that $\mathcal{A}_{K_x}^*$ is compact. Similarly, the stability properties of the signal generator \eqref{eq:automaton} and $w$ are asserted with respect to the set $[0,N_{0,u}]\times[0,T_{0,u}]\times\Sigma_u$ and $\mathcal{W}$, which are forward pre-invariant and trivially attractive since, by construction of the flow and jump set, no solution can start outside $[0,N_{0,u}]\times[0,T_{0,u}]\times\Sigma_u$ or $\mathcal{W}$ respectively.
\end{remark}
\vspace{0.1cm}
Next, we incorporate the dynamics of the plant into the stability analysis of the closed-loop system. In this case, we provide an upper bound for the gain of the controller, which guarantees exponential stability of the optimal point $z^*=(x^*,u_{w}^*)$ for the case when $\dot{w} = 0$.
\begin{proposition}\label{pro:st_ctrlattacked2}
Suppose that Assumptions \ref{as:cts_w}-\ref{as:attackModel} hold, and
let $\kappa_{u,1}$, and $\kappa_{u,2}$ satisfy the bound \eqref{timeratiocondition2}. Let the gain $\varepsilon$  satisfy:
\begin{align}\label{eq:eps_u}
    0<\varepsilon
    < \varepsilon^*:=\frac{\rho\underline{\lambda}(R)\min\{\mu,1\}}{2\ell_y \max\left\{\bar{M},1\right\}\norm{C}\norm{G}\norm{PA^{-1}B}(\rho\min\{\mu,1\} + 2\ell\max\left\{\ell,\bar{M}\right\}e^{\tau_0})},
\end{align}
where $\rho>0$, and $\tau_0 = \ln(\omega)N_{0,u} + 2T_{0,u}(\mu + \bar{M}\ell)$, and $P$ is defined as in Assumption \ref{as:stabilityPlant}. Then, the set 
\blue{
\begin{equation}\label{optimalset}
\mc A^* =\left\{(x,u,\tau_{u,1},\tau_{u,2},\sigma_u,w):
x=x^*, u=u_w^*, ~\tau_{u,1}\in[0,N_{0,u}],\tau_{u,2}\in [0,T_{0,u}],\sigma_u\in \Sigma_u, w \in \mathcal{W}\right\},
\end{equation}
}
is UGES for the HDS \eqref{eq:intercon_u} with $\dot{w}=0$.
\end{proposition}
The result of Proposition \ref{pro:st_ctrlattacked2} establishes a solution to Problem \ref{prob:controlObjective} under a constant $w$. Note that the result is actually uniform in $w$, i.e., the time ratio $\kappa_{u,2}$ and the gain $\varepsilon$ are independent of $w$. Indeed, as we will show below, they are also independent of $\norm{\dot{w}(t)}_t$, which is an important property needed to avoid vanishing gains and/or shrinking safety margins under highly oscillating exogenous signals.

\vspace{0.1cm}
Since by Lemma \ref{wellposedHDS} the HDS \eqref{eq:intercon_u} is well-posed, the following additional robustness result can be asserted for the closed-loop system under attacks \emph{and} small disturbances $e$ acting on the states and dynamics. In general, such type of disturbances is unavoidable in practical applications due to measurements noise, numerical approximations, etc.
\begin{lemma}\label{robustness_result}
Consider the HDS \eqref{eq:intercon_u} with $\dot{w}=0$, and suppose that the conditions of Propositions \ref{pro:kappa_u2} and \ref{pro:st_ctrlattacked2} hold for the dwell-time parameter $\kappa_{u,1}$, the time-ratio parameter $\kappa_{u,2}$ and the gain $\varepsilon>0$. Then, there exists $c_1,c_2>0$ such that for each compact set of initial conditions $K\subset\mathbb{R}^{n}\times\mathbb{R}^m$, and for each $\delta>0$, there exists $\bar{e}>0$ such that for any measurable disturbance $e$ with $\sup_{t+j\geq0}|e(t,j)|\leq \bar{e}$ and every initial condition $\vartheta(0,0)\in K\times[0,N_{0,u}]\times[0,T_{0,u}]\times\Sigma \times \mathcal{W}$, the trajectories of the perturbed HDS
\begin{subequations}
\begin{align}&\vartheta + e\in \mathcal{C},~~~~~\dot{\vartheta}\in F_{\sigma_u}(\vartheta+e)+e,\\
    &\vartheta+e\in \mathcal{D},~~~\vartheta^+\in G_{\sigma_u}(\vartheta+e)+e,
\end{align}
\end{subequations}
satisfy the bound
\begin{equation}
|z(t,j)|_{\mathcal{A}^*}\leq c_1|z(0,0)|_{\mathcal{A}^*}\exp(-c_2(t+j))+\delta,~~~\forall~(t,j)\in\text{dom}(\vartheta),
\end{equation}
\end{lemma}
The robustness result of Lemma \ref{robustness_result} is not trivial. Indeed, for general hybrid dynamical systems, it is difficult to guarantee closeness of solutions between nominal and perturbed dynamics, even when the perturbations are arbitrarily small. On the other hand, when the HDS satisfies the regularity conditions of Lemma \ref{wellposedHDS}, stability properties turn out to be robust to small disturbances \cite[Ch.7]{Goebel:12}. This idea is at the core of our modeling framework \eqref{eq:intercon_u}, which subsumes the closed-loop system under attack as a well-posed hybrid dynamical system with sufficiently slow control dynamics.

\vspace{0.1cm}
By leveraging the results of Propositions \ref{pro:kappa_u2} and \ref{pro:st_ctrlattacked2}, we can now state the first \emph{main} result of this paper, which establishes uniform global exponential input-to-state stability (E-ISS) for the closed-loop system under attacks and time-varying $w$. In particular, under the time ratio established in Proposition \ref{pro:kappa_u2}, and the maximal gain $\varepsilon^*$ established in Proposition \ref{pro:st_ctrlattacked2}, we provide an explicit characterization of the ISS gain $\gamma$ that maps $\norm{\dot{w}(t)}_t$ to the radius of the residual set where the trajectories converge.
\begin{theorem}\label{thm:iss_u}
Let $\mathcal{A}^*$ be given by \eqref{optimalset}. Suppose that Assumptions \ref{as:cts_w}-\ref{as:attackModel} hold, and that conditions \eqref{timeratiocondition2} and \eqref{eq:eps_u} are also satisfied for $(\kappa_{u,1},\kappa_{u,2},\varepsilon)$. Then, the HDS \eqref{eq:intercon_u} renders E-ISS the set $\mc A^*$ with respect to $\dot{w}(t)$, with linear asymptotic gain $\gamma$ given by
\begin{align}\label{eq:gamma_u}
\gamma\left(\norm{\dot{w}(t)}_t\right):=\frac{1}{\varepsilon}\left[ \left(\frac{\max\left\{\theta \bar{\lambda}(P), (1-\theta)e^{\tau_0}\max\left\{\ell^2/2\mu,1/2\right\}\right\}}{\min\left\{\theta \underline{\lambda}(P),(1-\theta)\min\{\mu/2,1/2\} \right\}}\right)^{1/2} \frac{\norm{r}}{\underline{\lambda}(\Xi)k}\right]\norm{\dot{w}(t)}_t,
\end{align}
where $k \in (0,1)$ , $\Xi$ is a positive definite matrix, and  $r: = \left(2\theta \norm{PA^{-1}E},(1-\theta)e^{\tau_0}\max\left\{\ell_y\norm{H}\norm{G},U_{\bar{w}}\right\}\right)$ with
\begin{align*}
    \theta:= \frac{\ell_y \norm{G}\norm{C}\max\left\{\ell,\bar{M}\right\}}{\ell_y \norm{G}\norm{C}\max\left\{\ell,\bar{M}\right\} + 2\ell e^{-\tau_0}\norm{PA^{-1}B}\max\left\{\bar{M},1\right\}}.
\end{align*}
\end{theorem}
\begin{remark}
The gain $\gamma$ in \eqref{eq:gamma_u} reflects explicitly the stability of the plant through the eigenvalues of the matrix $P$ set in Assumption \ref{as:stabilityPlant}, and the qualitative behavior of the cost function given by the Lipschitz and $PL$ constant $\ell$ and $\mu$, respectively. Additionally, this gain is weighted by $1/\varepsilon$, which implies that a slower controller with respect to the dynamics of the plant leads to larger tracking errors. This behavior is expected since slowing down the controller makes it more difficult to track time-varying optimal points.
\end{remark}
\begin{remark}
When $\ell>>1$ and $\mu<<1$ (i.e., the cost function $f$ is ill-conditioned) the asymptotic gain simplifies to
\begin{align}\label{gain_theorem22}
\gamma\left(\norm{\dot{w}(t)}_t\right):=\frac{1}{\varepsilon}\Bigg[ \text{cond}(f)e^{\tau_0/2} \frac{\norm{r}}{\underline{\lambda}(\Xi)k}\Bigg]\norm{\dot{w}(t)}_t,~~~~~~\text{cond}(f)=\frac{\ell}{\mu},
\end{align}
which shows the effect that the condition number of the cost $f(\cdot)$ has on the residual set where the trajectories converge. Note that a cost function $f$ that is ill-conditioned can be identified by examining the eccentricity of its sub-level sets.
\end{remark}
\begin{remark}
In \eqref{eq:gamma_u}, the constant $\tau_0$ comes from the Lyapunov-based analysis used to study the HDS \eqref{eq:intercon_u}. This is an upper bound for the weighted time that the closed-loop system \eqref{eq:interconnection} can flow under attacks. A detailed characterization of these constants is presented in Section \ref{section_analysis}. 
\end{remark}
\subsection{Attacks to the Static Feedback}

\vspace{0.1cm}
In the previous section, we focused on attacks acting on the dynamic controller. In this section, we now turn our attention to study the effect of persistent attacks acting on the static controller \eqref{eq:attackModelStatic}. To study this scenario, we consider the following set-valued HDS with state $\vartheta:=(x,u,\tau_{v,1},\tau_{v,2},\sigma_v,w)$, which corresponds to the interconnection of the plant dynamics under attacks \eqref{eq:attackModelStatic}, the nominal dynamic gradient-based controller \eqref{eq:attackModelDynamic}, and the signal generator \eqref{eq:automaton}:
\blue{
\begin{subequations}\label{eq:intercon_v}
\begin{align}
\vartheta \in \mathcal{C}&:=  \R^n \times \R^m \times [0,N_{0,v}] \times [0,T_{0,v}] \times \Sigma_v \times \mathcal{W}, \\
\begin{pmatrix} \dot{x}\\\dot{u}\\\dot{\tau}_{v,1}\\\dot{\tau}_{v,2} \\ \dot{\sigma}_v \\\dot{w} \end{pmatrix} &\in F_{\sigma_v}(\vartheta) := \begin{pmatrix}
    A_{\sigma_v} x + B_{\sigma_v} u + E_{\sigma_v} w\\
    - \varepsilon(\nabla f_u(u) + G^\top \nabla f_y(y))\\
    [0,\kappa_{v,1}]\\
    [0,\kappa_{v,2}] - I_{\sigma_{v,a}} ({\sigma_v})\\
    0\\ \Pi(w)\end{pmatrix}, \quad y = Cx,\label{flowresult1}\\
\vartheta \in \mathcal{D}& :=  \R^n \times \R^m \times [1,N_{0,v}] \times [0,T_{0,v}] \times \Sigma_v \times \mathcal{W},\\
\begin{pmatrix} x^+ \\ u^+\\ \tau_{v,1}^+\\\tau_{v,2}^+ \\ {\sigma_v}^+ \\w^+ \end{pmatrix} &\in G_{\sigma_v}(\vartheta) := \begin{pmatrix}
    x\\
    u\\
    \tau_{v,1}-1\\
    \tau_{v,2}\\
    \Sigma_v \backslash \{\sigma_v\}\\
    w
    \end{pmatrix},
\end{align}
\end{subequations}
}
where $T_{0,v}\geq 0$, $\kappa_{v,1}>0$, $\kappa_{v,2} \in (0,1)$ is a time-ratio parameter that describes an admissible attack, \blue{and $\Pi(w)$ is given by Assumption \ref{as:cts_w}. Note that when $\sigma_v = s$, we have $A_{s}= A$, hence, Assumption \ref{as:stabilityPlant} holds for $A_s$. Moreover, to simplify notation we denote $A_{a_i} : = F + N L_{a_i} K C$ for any $a_i \in \Sigma_{v,a}$, thus we will refer to the attacks to the static controller also as attacks to the plant.} By construction, system \eqref{eq:intercon_v} also satisfies the properties of Lemma \ref{wellposedHDS}. 

For attacks acting on the plant, we introduce an additional technical assumption to guarantee that the switching system has a well-defined unique equilibrium point.
\begin{assumption}\label{as:commonequilibrium}
For any $w \in \mathcal{W}$  \red{and $\bar{u} \in \R^m$ there exists a unique $\bar{x} \in \R^n$ such that} $0 = A_{\sigma} \bar{x} + B_{\sigma} \bar{u} + E_{\sigma} w,~\forall \sigma \in \Sigma $.
\end{assumption}
\begin{remark}
The conditions of Assumption \ref{as:commonequilibrium} assert that the attacks to the plant modify the stability of the equilibria, but they do not alter the set of equilibria~\cite{Senejohnny2018TCN}. We note that this assumption is aligned with existing works in context~\cite{Senejohnny2018TCN,WANG2020AUT}, and it enables a tractable analysis via Lyapunov-based tools. Without this assumption in place, each attack could generate a different equilibrium point (if any), and to assert stability properties we will need to introduce additional restrictive assumptions on the frequency of the attacks, e.g., condition \eqref{eq:avgdwelltimeCondition} with a ``sufficiently small'' constant $\kappa_{v,1}$. The study of switched dynamical systems with multiple equilibria (without attacks) is addressed in \cite{alpcan2010multieq, Veer2020MultipleEq} to mention a few, and extending our framework to this case is the subject of future research efforts.
\end{remark}

In contrast to the analysis of the previous section, we now first consider the plant under attacks, operating with a constant controller command $u$ (i.e., $\dot{u}=0$) and a fixed \blue{ $w$ (i.e., $\dot{w}=0$)}. Thus, we study the stability properties of the equilibrium point $x^*(u,w):=A_{\sigma_v}^{-1}B_{\sigma_v}u+A_{\sigma_v}^{-1}E_{\sigma_v}w$. For this reduced system under attacks, we establish the following characterization of the parameters $(\kappa_{v,2},\kappa_{v,1})$ to guarantee uniform global exponential stability.

\begin{proposition}\label{pro:kappa_v2}
Suppose that Assumptions 
\ref{as:cts_w}-\ref{as:commonequilibrium} hold. Given a $\tilde{\beta}>0$, and the dwell-time parameter $\kappa_{v,1}$, and the time-ratio parameter $\kappa_{v,2}$ simultaneously satisfying: 
\begin{align}\label{timeratiocondition1}
\underline{\lambda}(R)>\bar{\lambda}(P)\ln(\omega)\kappa_{v,1} + \kappa_{v,2}\left(\underline{\lambda}(R)+\bar{\lambda}(P)\bar{\lambda}(\hat{R}_a)\right),~~~\text{with}~~~\omega:= \max\left\{\frac{\bar{\lambda}(P)}{\tilde{\beta}},\frac{\tilde{\beta}}{\underline{\lambda}(P)}\right\}.
\end{align}
where $P,R$ satisfy Assumption \ref{as:stabilityPlant}, and $\bar \lambda(\hat{R}_a)$ is a constant derived after Lemma \ref{lm:symR}. Then,
for each each compact set $K_u\subset\mathbb{R}^m$, the set 
\blue{
$$\mc{A}^*_{K_u} =\left\{(x,u,\tau_{v,1},\tau_{v,2},\sigma_v,w):x=x^*(u,w),~u\in K_u,~\tau_{v,1}\in[0,N_{0,v}],\tau_{v,2}\in [0,T_{0,v}],\sigma_v\in \Sigma_v,w \in \mathcal{W}\right\},$$ 
is UGES for the HDS \eqref{eq:intercon_v} with $\varepsilon=0$, $\dot{w}=0$, and restricted flow 
and jump sets  given by
\begin{align*}
\hspace{-0.5cm}
    \mathcal{C}_{K_u} =  \R^n \times (\R^m\cap K_u) \times [0,N_{0,v}] \times [0,T_{0,v}] \times \Sigma_v \times \mathcal{W},~\mathcal{D}_{K_u} =  \R^n \times (\R^m\cap K_u) \times [1,N_{0,v}] \times [0,T_{0,v}] \times \Sigma_v \times \mathcal{W}.
\end{align*}
}
\end{proposition}
\vspace{-0.5cm}
The result of Proposition \ref{pro:kappa_v2} establishes an explicit characterization of the admissible attack properties (c.f. Definition \ref{def:undetectableAttack}) needed to guarantee exponential stability in the system. For the internal feedback loop controller, these properties are related to the largest and smallest eigenvalues of the matrices $P$ and $\hat{R}_a$. 

\vspace{0.1cm}
By leveraging the result of Proposition \ref{pro:kappa_v2}, we can now state a stability and convergence result for the case where the controller is activated, and the state $w$ remains static. 
\begin{proposition}\label{pro:st_plantattacked}
Suppose that Assumptions \ref{as:cts_w}-\ref{as:commonequilibrium} hold, and
let $\kappa_{v,1}$, $\kappa_{v,2}$ satisfy the bounds in  Proposition 
\ref{pro:kappa_v2}. Let $\varepsilon$  satisfy:
 \begin{align}\label{eq:eps_v}
    0<\varepsilon < \varepsilon^*:=  \frac{\rho \min\{\underline{\lambda}(P),\tilde{\beta}\}\mu^2}{2\ell_y e^{\tau_0}\norm{C}\norm{G}\norm{A^{-1}B}\max\left\{\bar{\lambda}(P),\tilde{\beta}\right\}(\mu^2 + \ell^2)},
\end{align}
where $\rho>0$, $\tau_0 = \ln(\omega)N_{0,v} + T_{0,v}\left(\frac{\underline{\lambda}(R)}{\bar{\lambda}(P)} + \bar{\lambda}(\hat{R}_a)\right)$, and $P$
is defined in Assumption \ref{as:stabilityPlant}. Then, the set \blue{ 
$$\mc A^* =\left\{(x,u,\tau_{v,1},\tau_{v,2},\sigma_v,w):
x=x^*(u^*,w), u=u^*, ~\tau_{v,1}\in[0,N_{0,v}],\tau_{v,2}\in [0,T_{0,v}],\sigma_v\in \Sigma_v,w \in \mathcal{W}\right\},$$} 
is UGES for the HDS \eqref{eq:intercon_v} with $\dot{w}=0$.
\end{proposition}

\begin{remark}
The result of Proposition \ref{pro:st_plantattacked} establishes a sufficient condition on the bounds of the gain of the controller, the dwell-time parameter $\kappa_{v,1}$, and the time-ratio parameter $\kappa_{v,2}$ to obtain exponential convergence to the solution of Problem \eqref{prob:controlObjective}. Since the HDS is well-posed, the robustness results of Lemma \ref{robustness_result} also hold for system \eqref{eq:intercon_v}. 
\end{remark}
By leveraging the results of Proposition \ref{pro:kappa_v2} and Proposition \ref{pro:st_plantattacked}, we now establish the second main result of this paper, which asserts exponential input-to-state stability of the closed-loop system under attacks, with respect to the time-variation of $w$. In conjunction with Theorem \ref{thm:iss_u}, this result provides a complete answer to Problem 1.
\begin{theorem}\label{thm:iss_v}
Suppose that Assumptions \ref{as:cts_w}-\ref{as:commonequilibrium} hold, and that conditions \eqref{timeratiocondition1} and \eqref{eq:eps_v} also hold. Then, the HDS \eqref{eq:intercon_v} renders the set $\mc A^*$ E-ISS with respect to $\dot{w}(t)$, with linear asymptotic gain $\gamma$ given by
\begin{align}\label{gain_theorem2}
\gamma\left(\norm{\dot{w}(t)}_t\right):=\frac{1}{\varepsilon}\Bigg[ \left(\frac{\max\left\{\theta e^{\tau_0}\max\left\{\bar{\lambda}(P),\tilde{\beta}\right\},(1-\theta)\ell^2/(2\mu)\right\}}{\min\left\{\theta \min\left\{\underline{\lambda}(P),\tilde{\beta}\right\},(1-\theta)\mu/2\right\}}\right)^{1/2} \frac{\norm{r}}{\underline{\lambda}(\Xi)k}\Bigg]\norm{\dot{w}(t)}_t,
\end{align}
where $k\in (0,1)$, $\Xi$ is a positive definite matrix, and $r: = \left(2e^{\tau_0}\theta \max\{\bar{\lambda}(P),\tilde{\beta}\}\norm{A^{-1}E}, \ell_y (1-\theta)\norm{H}\norm{G}\right)$ with

\begin{align*}
    \theta:= \frac{\ell_y \norm{G}\norm{C}}{\ell_y \norm{G}\norm{C} + 2e^{\tau_0}\max\left\{\bar{\lambda}(P),\tilde{\beta}\right\}\norm{A^{-1}B}}.
\end{align*}
\end{theorem}
\begin{remark}
As in Theorem \ref{thm:iss_u}, the asymptotic gain \eqref{gain_theorem2} relates the magnitude of $\dot{w}$ to the size of the residual set where the trajectories converge. Note that in this case the coefficient $e^{\tau_0}$ appears next to the the values and gains related to the plant. Therefore, when $\bar{\lambda}(P)>>1$ and $\underline{\lambda}(P)<<1$, the asymptotic gain simplifies to
\begin{align}\label{eq:gamma_u2}
\gamma\left(\norm{\dot{w}(t)}_t\right):=\frac{1}{\varepsilon}\Bigg[ \sqrt{\text{cond}(P)e^{\tau_0}} \frac{\norm{r}}{\underline{\lambda}(\Xi)k}\Bigg]\norm{\dot{w}(t)}_t,~~~~~~\text{cond}(P):=\frac{\bar{\lambda}(P)}{\underline{\lambda}(P)},
\end{align}
which shows the effect of the condition number of the Lyapunov matrix $P$ on the residual set where the trajectories converge.
\end{remark}

\begin{figure}
    \centering
    \includegraphics[width=0.95\columnwidth]{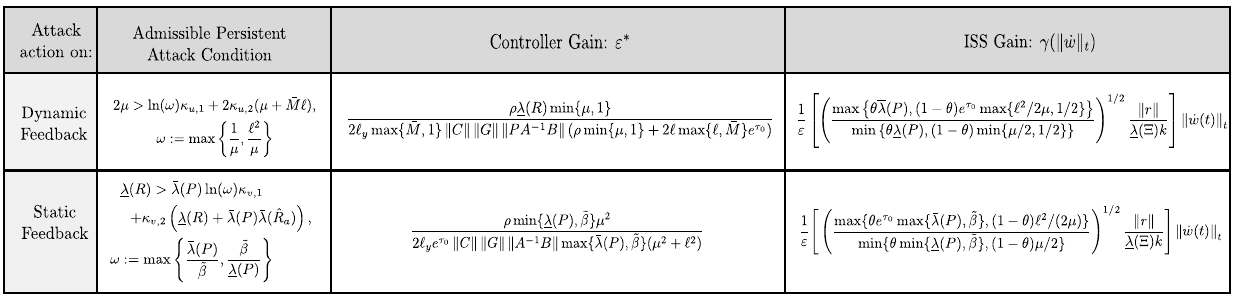}
    \caption{Summary of the results of Propositions \ref{pro:kappa_u2}-\ref{pro:st_plantattacked}, and Theorems \ref{thm:iss_u}-\ref{thm:iss_v}.}
    \label{fig:summaryTable}
\end{figure}

We summarize in Table \ref{fig:summaryTable} the main results of this section. Each row indicates which subsystem is under attack, and the columns summarize the upper bounds on the time ratios and gains, and also the form of the ISS gain.  

\section{Analysis}
\label{section_analysis}
In this section, we present the analysis and the proofs of Theorems~\ref{thm:iss_u} and~\ref{thm:iss_v}. In both cases, our analysis leverages Lyapunov theory for switched systems with unstable subsystems. The structure of the proofs will follow four main steps:
\begin{enumerate}
    \item For each of the feedback loops under attack, we will introduce an auxiliary Lyapunov-like function for which a strong exponential decrease can be asserted whenever the feedback loop operates in the stable mode, and the dynamics of the other loop are neglected.
    \item When the feedback loop is under attacks, we incorporate the switching behavior into the system by extending the auxiliary Lyapunov-like function with a state that depends on the hybrid signal generator.
    \item We use the extended auxiliary functions constructed in the previous step, to construct a new Lyapunov function for the complete hybrid system, and we characterize the upper bound $\varepsilon^*$ on the gain of the controller needed to guarantee exponential decrease of the function during flows outside of the compact set of interest. We also show that the Lyapunov function does not increase during jumps. 
    \item The previous argument, in addition to the quadratic upper and lower bounds on the Lyapunov function, and the average dwell-time constraint, allows us to establish uniform global exponential stability and/or uniform global exponential ISS.
\end{enumerate}
%

\subsection{Analysis of System with Dynamic Feedback Controller under Attacks}

\vspace{0.1cm}
Consider the closed-loop system with attacks on the gradient-flow controller modeled by the HDS \eqref{eq:intercon_u}. Define the affine map $(u,w) \mapsto \bar{x}(u,w)$ as 
$\bar{x}(u,w) := - A^\inv B u - A^\inv E w$, and consider the change of variables $x_e := x - \bar{x}$. In the new variables, the dynamics are given by:
\begin{align*}
 &\dot{x}_e = A x_e +  A^{-1}B\dot{u} + A^{-1}E\dot{w}, \quad\quad \dot{u} =- \varepsilon M_{\sigma_u}\varphi(u,w,x_e)= - \varepsilon M_{\sigma_u} \left(\nabla_u f_u(u) + G^\top \nabla_y f_y(C x_e + Gu + Hw)\right),\\
 &G:=-C A^\inv B, \quad\quad H:= -C A^\inv E,
\end{align*}
where we note that $-\varphi(u,w,0) = -\nabla_u f(u,w) = -\nabla_u f_u(u) - G^\top \nabla_y f_y( Gu + Hw)$. First, we investigate the stability properties of the stable mode of the system by leveraging the stability properties of the gradient-flow. Consider the following auxiliary function:
\begin{align*}
V_s(u,w) &= f(u,w) - f(u^*_w,w).
\end{align*}
For $ \vartheta \in \mathcal{D}$, during the jumps it holds that $V_s(u^+,w^+) = V_s(u,w)$. Moreover, when $\vartheta \in \mathcal{C}$, during flows we have:
\begin{align*}
\dot V_s(u,w) &= \nabla_u f(u,w)^\top \dot{u} + (H \Pi(w))^\top\nabla_y f_y(Gu+Hw) - \nabla_{u} f(u^*_w,w)^\top \dot{u}^* - (H \Pi(w))^\top\nabla_y f_y(Gu^*_w+Hw)\\
    &= \nabla_u f(u,w)^\top \dot{u} + (H \Pi(w))^\top\nabla_y f_y(Gu+Hw)  - (H \Pi(w))^\top\nabla_y f_y(Gu^*_w+Hw)\\
    &= - \varepsilon \varphi(u,w,0)^\top (\varphi(u,w,0) + \varphi(u,w,x_e) - \varphi(u,w,0)) + \ell_y \norm{H}\norm{G}\norm{u - u^*_w}\norm{\Pi(w)}\\
    &\leq - \varepsilon \norm{\nabla_u f(u,w)}^2
        + \varepsilon \ell_y \norm{C} \norm{G} \norm{\nabla_u f(u,w)} \norm{x_e}+ \ell_y \norm{H}\norm{G}\norm{u - u^*_w}\norm{\Pi(w)}\\
    &\leq - 2 \mu \varepsilon V_s(u,w) + \varepsilon\ell_y \norm{C} \norm{G} \norm{\nabla_u f(u,w)} \norm{x_e}+ \ell_y \norm{H}\norm{G}\norm{u - u^*_w}\norm{\Pi(w)}.
\end{align*}

Next, to study the stability properties of the controller under switching signals with unstable modes, we consider the following auxiliary function:
\begin{align*}
  V_a(u,w) &= \frac{1}{2}\norm{u - u^*_w}^2.
\end{align*}
When $\vartheta \in \mathcal{D}$, during jumps it holds that $V_a(u^+,w^+) = V_a(u,w)$. For $\vartheta \in \mathcal{C}$, during flows we have
\begin{align*}
\dot V_a(u,w) &= (u - u^*_w)^\top(\dot{u} - \dot{u}^*_w)\\
&= -(u - u^*_w)^\top\left(\varepsilon M_{\sigma_u}\varphi(u,w,x_e)\right) - (u - u^*_w)^\top(\nabla_w h \Pi(w))\\
&\leq \varepsilon \norm{u - u^*_w}\norm{M_{\sigma_u}}\norm{\varphi(u,w,x_e)} + \norm{u - u^*_w}\norm{\nabla_w h}\norm{\Pi(w)}\\
&\leq \varepsilon \norm{u - u^*_w}\norm{M_{\sigma_u}}\norm{\varphi(u,w,0)-\varphi(u,w,0)+\varphi(u,w,x_e)} + \norm{u - u^*_w}\norm{\nabla_w h}\norm{\Pi(w)}\\
&\leq \varepsilon \norm{u - u^*_w}\norm{M_{\sigma_u}}\left(\ell\norm{u-u^*_w} + \ell_y\norm{C}\norm{G}\norm{x_e}\right) + \norm{u - u^*_w}\norm{\nabla_w h}\norm{\Pi(w)}\\
&\leq \varepsilon \ell \bar{M}\norm{u - u^*_w}^2 + \varepsilon \bar{M}\ell_y\norm{C}\norm{G}\norm{u - u^*_w}\norm{x_e} + \norm{u - u^*_w}\norm{\nabla_w h}\norm{\Pi(w)}\\
&= 2\varepsilon \ell \bar{M}V_a(u,w) + \varepsilon \bar{M}\ell_y\norm{C}\norm{G}\norm{u - u^*_w}\norm{x_e} + U_{\bar{w}}\norm{u - u^*_w} \norm{\Pi(w)}.
\end{align*}
\blue{
Next, we consider the Lyapunov function: 
\begin{equation}
\begin{aligned}
\hspace{-0.3cm}
    \bar{V} (u,w,\tau_{u,1},\tau_{u,2},\sigma_u):= V_{\sigma_u}(u,w) e^\tau,~~ \sigma_u \in \{s,a\},~~\tau:= \ln(\omega)\tau_{u,1} + \tau_{u,2} (\rho_s + \rho_a),~~\rho_s := 2 \mu,~\rho_a :=  2\bar{M}\ell,
\end{aligned}
\end{equation}
where we recall that $\mu$ and $\ell$ are the PL and Lipschitz constants, respectively, $\bar{M}$ is given by Assumption \ref{as:attackModel}, and $\omega =  \max\left\{1/\mu,\ell^2/\mu \right\}$. Now, during jumps the function $\bar{V}(u,w,\tau_{u,1},\tau_{u,2},\sigma_u)$ satisfies
\begin{align*}
\bar{V}(u^+,w^+,\tau_{u,1}^+, \tau_{u,2}^+, \sigma_u^+) = V_{\sigma_u^+}e^{\ln(\omega)(\tau_{u,1}-1) + \tau_{u,2} (\rho_s + \rho_a)}
    \leq \max_{\sigma_u^+ \in \{s,a\}} V_{\sigma_u^+}e^{-\ln(\omega) + \tau} \leq \omega V_{\sigma_u}e^{-\ln(\omega)}e^{\tau} =  \bar{V},
\end{align*}
hence, we conclude that $\bar{V}$ does not increase during jumps.} 
Let $\tau_0 : = \ln(\omega)N_{0,u} + T_{0,u} (\rho_s + \rho_a)$. To analyze the flows, we note that for $\bar{V} = V_{s}e^\tau$ the time derivative satisfies:
\begin{align*}
\hspace{-0.9cm}\dot{\bar{V}}(u,w,\tau_{u,1},\tau_{u,2},\sigma_u)
&= \dot{V}_s(u,w)e^\tau + V_s(u,w)e^\tau \dot{\tau}\\
&= \dot{V}_s(u,w)e^\tau +  V_s(u,w)e^\tau(\ln(\omega)\dot{\tau}_{u,1} + (\rho_s + \rho_a)\dot{\tau}_{u,2}) \\
&\leq \dot{V}_s(u,w)e^\tau + \varepsilon V(u,w)_1e^\tau(\ln(\omega)\kappa_{u,1} + (\rho_s + \rho_a)\kappa_{u,2}) \\
&\leq - 2 \mu \varepsilon V_s(u,w)e^\tau + \varepsilon\ell_y e^\tau\norm{C} \norm{G} \norm{\nabla_u f(u,w)} \norm{x_e}\\
    &\quad + \ell_ye^\tau \norm{H}\norm{G}\norm{u - u^*_w}\norm{\Pi(w)} +  \varepsilon V_s(u,w)e^\tau (\ln(\omega)\kappa_{u,1} + (\rho_s + \rho_a)\kappa_{u,2})  \\
 &= -\varepsilon\rho  \bar{V}(u,w,\tau_{u,1},\tau_{u,2},\sigma_u) + \varepsilon\ell_y e^\tau\norm{C} \norm{G} \norm{\nabla_u f(u,w)} \norm{x_e} + \ell_y e^\tau \norm{H}\norm{G}\norm{u - u^*_w}\norm{\Pi(w)}\\
  &\leq -\varepsilon\rho  \bar{V}(u,w,\tau_{u,1},\tau_{u,2},\sigma_u) + \varepsilon\ell \ell_y e^{\tau_0}\norm{C} \norm{G} \norm{u - u^*_w} \norm{x_e} + \ell_y e^{\tau_0} \norm{H}\norm{G}\norm{u - u^*_w}\norm{\Pi(w)}\\
  &\leq -\varepsilon\rho  \bar{V}(u,w,\tau_{u,1},\tau_{u,2},\sigma_u) + \varepsilon\ell_y e^{\tau_0}\max\left\{\ell,\bar{M}\right\} \norm{C} \norm{G} \norm{u - u^*_w} \norm{x_e}\\
   &\quad + e^{\tau_0} \max\left\{\ell_y \norm{H}\norm{G}, U_{\bar{w}}\right\}\norm{u - u^*_w}\norm{\Pi(w)},
\end{align*}
where $\rho := \rho_s - \kappa_{u,2} (\rho_s + \rho_a) - \ln(\omega)\kappa_{u,1} >0$, for $\omega \geq 1$. Similarly, when $\bar{V} = V_a e^\tau$ we have that:
\begin{align*}
\hspace{-0.9cm}\dot{\bar{V}}(u,w,\tau_{u,1},\tau_{u,2},\sigma_u)
&= \dot{V}_a(u,w)e^\tau + V_a(u,w)e^\tau \dot{\tau}\nonumber\\
&= \dot{V}_a(u,w)e^\tau +  V_a(u,w)e^\tau(\ln(\omega)\dot{\tau}_{u,1} + (\rho_s + \rho_a)\dot{\tau}_{u,2}) \nonumber\\
&\leq \dot{V}_a(u,w)e^\tau + \varepsilon V_a(u,w)e^\tau(\ln(\omega)\kappa_{u,1} + (\rho_s + \rho_a)(\kappa_{u,2}-1)) \nonumber\\
&\leq 2\varepsilon \ell \bar{M}V_a(u,w)e^\tau + \varepsilon \bar{M}\ell_ye^\tau\norm{C}\norm{G}\norm{u - u^*_w}\norm{x_e} + U_{\bar{w}}\norm{u - u^*_w} \norm{\Pi(w)}e^\tau \\
    &\quad + \varepsilon V_a(u,w)e^\tau(\ln(\omega)\kappa_{u,1} + (\rho_s + \rho_a)(\kappa_{u,2}-1))\\
& =  -\varepsilon\rho  \bar{V}(u,w,\tau_{u,1},\tau_{u,2},\sigma_u) + \varepsilon \bar{M}\ell_y e^{\tau_0}\norm{C}\norm{G}\norm{u - u^*_w}\norm{x_e} + U_{\bar{w}}\norm{u - u^*_w} \norm{\Pi(w)}e^{\tau_0}\\
 &\leq -\varepsilon\rho  \bar{V}(u,w,\tau_{u,1},\tau_{u,2},\sigma_u) + \varepsilon\ell_y e^{\tau_0}\max\left\{\ell,\bar{M}\right\} \norm{C} \norm{G} \norm{u - u^*_w} \norm{x_e}\\
   &\quad + e^{\tau_0} \max\left\{\ell_y \norm{H}\norm{G}, U_{\bar{w}}\right\}\norm{u - u^*_w}\norm{\Pi(w)},
\end{align*}
Therefore, for any mode we can upper-bound $\dot{\bar{V}}$ as follows:
\begin{align*}
\hspace{-0.9cm}\dot{\bar{V}}(u,w,\tau_{u,1},\tau_{u,2},\sigma_u)
 &\leq -\varepsilon\frac{\rho}{2}\min\{\mu,1\}\norm{u - u^*_w}^2 + \varepsilon\ell_y e^{\tau_0}\max\left\{\ell,\bar{M}\right\} \norm{C} \norm{G} \norm{u - u^*_w} \norm{x_e}\\
   &\quad + e^{\tau_0} \max \left\{\ell_y \norm{H}\norm{G}, U_{\bar{w}}\right\}\norm{u - u^*_w}\norm{\Pi(w)}.
\end{align*}
We now analyze the interconnection of the plant dynamics and the switched controller. To do this, we first consider the auxiliary function $W(x_e) = x_e^\top P x_e$, where $P \succ 0$ which satisfies $A^\top P + P A \preceq -R$ with $R \succ 0$ due to Assumption \ref{as:stabilityPlant}. When $\vartheta \in \mathcal{D}$, during jumps we have $W(x_e^+) = W(x_e)$. Moreover, when $\vartheta \in \mathcal{C}$, the time derivative of $W$ along the trajectories of the plant dynamics satisfies the following inequalities:
\begin{align*}
\hspace{-1cm}
\dot{W}(x_e) &= \dot{x}_e^\top Px_e + x_e^\top P \dot{x}_e\nonumber\\
&= 2x_e^\top P(A x_e +  A^{-1}B\dot{u} + A^{-1}E\Pi(w))\\
& = x_e^\top (A^\top P + PA)x_e + 2x_e^\top\left( P A^{-1}B \dot{u}\right) + 2x_e^\top\left( P A^{-1}E \Pi(w)\right)\nonumber\\
&\leq -x_e^\top R x_e - 2\varepsilon x_e^\top\left( P A^{-1}B M_{\sigma_u}\varphi(u,w,x_e)\right)+ 2x_e^\top\left( P A^{-1}E \Pi(w)\right)\nonumber\\
&\leq -\underline{\lambda}(R)\norm{x_e}^2 + 2 \varepsilon \norm{PA^{-1}B}\norm{M_{\sigma_u}}\norm{\varphi(u,w,x_e)}\norm{x_e} +2\norm{P A^{-1}E}\norm{x_e}\norm{\Pi(w)}\nonumber \\
& \leq -\underline{\lambda}(R)\norm{x_e}^2 + 2 \varepsilon \norm{PA^{-1}B}\norm{M_{\sigma_u}}\norm{x_e}\norm{\varphi(u,w,0)+\varphi(u,w,x_e)-\varphi(u,w,0)}\\ 
&\quad +2\norm{P A^{-1}E}\norm{x_e}\norm{\Pi(w)} \nonumber \\
& \leq -\underline{\lambda}(R)\norm{x_e}^2 + 2 \varepsilon \norm{PA^{-1}B}\norm{M_{\sigma_u}}\norm{x_e}\left(\norm{\varphi(u,w,0)}+\norm{\varphi(u,w,x_e)-\varphi(u,w,0)}\right)\\
&\quad +2\norm{P A^{-1}E}\norm{x_e}\norm{\Pi(w)} \nonumber \\
& \leq -\underline{\lambda}(R)\norm{x_e}^2 + 2 \varepsilon \norm{PA^{-1}B}\norm{M_{\sigma_u}}\norm{x_e}\left(\norm{\nabla_u f(u,w)}+\ell_y\norm{C}\norm{G}\norm{x_e}\right) +2\norm{P A^{-1}E}\norm{x_e}\norm{\Pi(w)} \nonumber \\
&= -\underline{\lambda}(R)\norm{x_e}^2 +2\varepsilon \ell_y  \norm{PA^{-1}B}\norm{M_{\sigma_u}}\norm{C}\norm{G}\norm{x_e}^2 + 2 \varepsilon \norm{PA^{-1}B}\norm{M_{\sigma_u}}\norm{x_e}\norm{\nabla_u f(u,w)}\\
&\quad +2\norm{P A^{-1}E}\norm{x_e}\norm{\Pi(w)} \nonumber\\
&\leq -\underline{\lambda}(R)\norm{x_e}^2 +2\varepsilon \ell_y \max\left\{\bar{M},1\right\} \norm{PA^{-1}B}\norm{C}\norm{G}\norm{x_e}^2 + 2 \varepsilon \max\left\{\bar{M},1\right\} \norm{PA^{-1}B}\norm{x_e}\norm{\nabla_u f(u,w)}\\
&\quad + 2\norm{P A^{-1}E}\norm{x_e}\norm{\Pi(w)}. \nonumber\\
&\leq -\underline{\lambda}(R)\norm{x_e}^2 +2\varepsilon \ell_y \max\left\{\bar{M},1\right\} \norm{PA^{-1}B}\norm{C}\norm{G}\norm{x_e}^2 + 2 \varepsilon \ell\max\left\{\bar{M},1\right\} \norm{PA^{-1}B}\norm{x_e}\norm{u - u^*_w}\\
&\quad + 2\norm{P A^{-1}E}\norm{x_e}\norm{\Pi(w)}. \nonumber
\end{align*}
Next, we consider the following Lyapunov function for the complete HDS \eqref{eq:intercon_u}: 
\blue{
\begin{equation}
U(\vartheta) = (1-\theta) \bar V(u,w,\tau_{u,1}, \tau_{u,2},\sigma_u) + \theta W(x_e),~~~~\theta\in(0,1).
\end{equation}
}
This function is bounded from below and above as follows:
\begin{equation}\label{eq:quadbounds_Uu}
\min\left\{\theta \underline{\lambda}(P),(1-\theta)\frac{1}{2}\min\{\mu,1\} \right\}|\vartheta|^2_{\mathcal{A}^*}\leq U(\vartheta)\leq \max\left\{\theta \bar{\lambda}(P), (1-\theta)\frac{1}{2}\max\left\{\ell^2/\mu,1\right\}e^{\tau_0}\right\}|\vartheta|^2_{\mathcal{A}^*},
\end{equation}
$\forall~\vartheta \in \mathcal{C}\cup \mathcal{D}$, where we used the fact that $\tau_{u,1}\in[0,N_{0,u}]$, $\tau_{u,2}\in[0,T_{0,u}]$, $\sigma_u\in \Sigma_u$, and $w \in \mathcal{W}$ at all times. Moreover, for any $\vartheta \in \mathcal{D}$, it follows that \blue{$U(\vartheta^+) \leq (1-\theta)\bar{V} + \theta W = U(\vartheta)$}, hence, during jumps, the Lyapunov function does not increase. For any $\vartheta \in \mathcal{C}$, the time derivative of $\dot U$  can be upper bounded as follows: 
\begin{align*}
\dot{U}(\vartheta) &= -\theta\underline{\lambda}(R)\norm{x_e}^2 +2\varepsilon \theta\ell_y \max\left\{\bar{M},1\right\} \norm{PA^{-1}B}\norm{C}\norm{G}\norm{x_e}^2\\
    &\quad + 2 \varepsilon \theta\ell\max\left\{\bar{M},1\right\} \norm{PA^{-1}B}\norm{x_e}\norm{u - u^*_w} + 2\theta\norm{P A^{-1}E}\norm{x_e}\norm{\Pi(w)}\\
    &\quad -\varepsilon(1-\theta)\frac{\rho}{2}\min\{\mu,1\}\norm{u - u^*_w}^2 + \varepsilon\ell_y e^{\tau_0}(1-\theta)\max\left\{\ell,\bar{M}\right\} \norm{C} \norm{G} \norm{u - u^*_w} \norm{x_e}\\
    &\quad + e^{\tau_0} (1-\theta)\max\left\{\ell_y \norm{H}\norm{G}, U_{\bar{w}}\right\}\norm{u - u^*_w}\norm{\Pi(w)}. 
\end{align*}
Let $\xi : = \left(\norm{x_e}, \norm{u - u^*_w}\right)$; using these definitions, we obtain: 
\begin{align}\label{eq:ISS_u1} 
    \dot{U}(\vartheta) \leq -\varepsilon \xi^\top \Xi \xi + r^\top \xi \norm{\Pi(w)},
\end{align}
where $r: = \left(2\theta \norm{PA^{-1}E}, e^{\tau_0} (1-\theta)\max\left\{\ell_y \norm{H}\norm{G}, U_{\bar{w}}\right\}\right)$, and $\Xi$ is a symmetric matrix of the form 
\begin{align}
    \Xi = \begin{bmatrix}
        \theta\left(\frac{\alpha}{\varepsilon} - \beta\right) & -\frac{1}{2}\left((1-\theta)\delta + \theta \chi\right)\\
        -\frac{1}{2}\left((1-\theta)\delta + \theta \chi\right) & (1-\theta)\gamma
    \end{bmatrix}.
\label{eq:matrix_Xi}
\end{align}
with $\alpha:= \underline{\lambda}(R)$, $\beta:= 2\ell_y \max\left\{\bar{M},1\right\}\norm{PA^{-1}B}\norm{C}\norm{G}$, $\delta := \ell_y
e^{\tau_0}\max\left\{\ell,\bar{M}\right\}\norm{G}\norm{C}$,

\noindent
$\chi:= 2\ell\max\left\{\bar{M},1\right\}\norm{PA^{-1}B}$, $\gamma:=  \frac{\rho}{2}\min\{\mu,1\}$, $\theta:=\frac{\delta}{\delta + \chi}$, and with  $\varepsilon$ satisfying \eqref{eq:eps_u},  hence, $\Xi$ is positive definite according to Lemma \ref{lm:posdefmatrix}.

The last step is to show that the quadratic term on the right-hand-side of~\eqref{eq:ISS_u1} dominates the linear term whenever $\xi$ is sufficiently large. It follows that for $k \in (0,1)$, \eqref{eq:ISS_u1} can be further upper bounded as follows:
\begin{align}\label{eq:ISS_u2}
    \dot{U}(\vartheta) &\leq -\varepsilon(1-k)\underline{\lambda}(\Xi)\norm{\xi}^2 -\varepsilon k\underline{\lambda}(\Xi)\norm{\xi}^2+ r^\top \xi \norm{\Pi(w)}_t\nonumber\\
    &\leq -\varepsilon(1-k)\underline{\lambda}(\Xi)\norm{\xi}^2 -\varepsilon k\underline{\lambda}(\Xi)\norm{\xi}^2+ \norm{r}\norm{\xi} \norm{\Pi(w)}_t\nonumber\\
    &\leq -\varepsilon(1-k)\underline{\lambda}(\Xi)\norm{\xi}^2=-\varepsilon (1-k)\underline{\lambda}(\Xi) |\vartheta|_{\mathcal{A}^*}^2,~~~\forall~|\vartheta|_{\mathcal{A}^*} \geq \frac{\norm{r}}{\varepsilon \underline{\lambda}(\Xi) k }\norm{\Pi(w)}_t.
\end{align}

By applying \cite[Thm. 3.1]{CAI200947}, exploiting the quadratic bounds in \eqref{eq:quadbounds_Uu}, and the average dwell-time condition, the result holds.
\hfill $\blacksquare$
%
\subsection{Analysis of System with Static Feedback Controller under Attacks}
In this section, we study the stability properties of system~\eqref{eq:intercon_v}, which models the interconnected system with attacks acting on the inner control loop. \blue{We introduce the affine mapping $(u,w) \mapsto \bar{x}(u,w)$ given by $\bar{x}(u,w) = -A_{\sigma_v}^{-1}B_{\sigma_v}u  - A_{\sigma_v}^{-1}E_{\sigma_v}w$, which, by Assumption \ref{as:commonequilibrium}, is independent of $\sigma_v \in \Sigma_v$. Using this mapping, we consider the change of variable $x_e : = x - \bar{x}$, which leads to the following dynamics:}
%
\begin{align*}
	&\dot{x}_e =  A_{\sigma_v} x_e +  A_{\sigma_v}^{-1}B_{\sigma_v} \dot{u} + A_{\sigma_v}^{-1}E_{\sigma_v} \dot{w},\qquad \dot{u} =  - \varepsilon\varphi(u,w,x_e): = - \varepsilon(\nabla f_u(u) + G_{\sigma_v}^\top \nabla f_y (Cx_e + G_{\sigma_v} u + H_{\sigma_v} w)),\\
	& G_{\sigma_v}:= -CA_{\sigma_v}^{-1}B_{\sigma_v}, \qquad H_{\sigma_v}:= -CA_{\sigma_v}^{-1}E_{\sigma_v},
\end{align*}
where $-\varphi(u,w,0)=-\nabla_u f(u,w) =  -\nabla f_u(u) - G_s^\top \nabla f_y(G_{\sigma_v} u +  H_{\sigma_v} w)$.

We first consider the stable mode, i.e., $\sigma_v=s$. For this mode, we consider the auxiliary function $W_s := x_e^\top P x_e$, where $P \succ 0$ satisfies the algebraic condition $A_s^\top P + P A_s = -R$ for a given matrix $R \succ 0$. When $\vartheta \in \mathcal{D}$, during the jumps we have that $W_s(x_e^+) = W_s(x_e)$. On the other hand, when $\vartheta \in \mathcal{C}$, the time derivative of $W_s$ satisfies:
\begin{align*}
    \dot{W}_s(x_e) &= \dot{x}_e^\top P x_e +  x_e^\top P \dot{x}_e\\
    & = \left(A_s x_e + A_s^{-1}B_s\dot{u} + A_s^{-1}E_s\Pi(w)\right)^\top P x_e + x_e^\top P \left(A_s x_e + A_s^{-1}B_s\dot{u} + A_s^{-1}E_s\Pi(w)\right)\\
    & = x_e^\top(A_s^\top P + P A_s)x_e + 2 x_e^\top P(A_s^{-1}B_s \dot{u}) + 2 x_e^\top P(A_s^{-1}E_s \Pi(w))\\
    & = -x_e^\top R x_e - 2\varepsilon x_e^\top P A_s^{-1}B_s(\varphi(u,w,x_e) + \varphi(u,w,0) - \varphi(u,w,0))+ 2 x_e^\top P(A_s^{-1}E_s \Pi(w))\\
    &\leq -x_e^\top R x_e + 2\varepsilon \norm{x_e} \norm{ P A_s^{-1}B_s}\norm{ \varphi(u,w,0) + (\varphi(u,w,x_e)  - \varphi(u,w,0))} + 2 x_e^\top P(A_s^{-1}E_s \Pi(w))\\
    &\leq -x_e^\top R x_e + 2\varepsilon \norm{x_e} \norm{ P A_s^{-1}B_s}\left(\norm{ \varphi(u,w,0)} + \norm{\varphi(u,w,x_e)  - \varphi(u,w,0)}\right)+ 2 x_e^\top P(A_s^{-1}E_s \Pi(w))\\
    & \leq -\underline{\lambda}(R)\norm{x_e}^2 + 2\varepsilon\norm{PA_s^{-1}B_s}\norm{x_e}\norm{\nabla_u f(u,w)} + 2\varepsilon\ell_y\norm{P A_s^{-1}B_s}\norm{C}\norm{G_s}\norm{x_e}^2\\
    &\quad + 2 \norm{x_e}\norm{ P A_s^{-1}E_s}\norm{ \Pi(w)}\\
    & \leq -\rho_s W_s(x_e) +  2\varepsilon\norm{PA_s^{-1}B_s}\norm{x_e}\norm{\nabla_u f(u,w)} + 2\varepsilon\ell_y\norm{P A_s^{-1}B_s}\norm{C}\norm{G_s}\norm{x_e}^2\\
    &\quad + 2 \norm{x_e}\norm{ P A_s^{-1}E_s}\norm{ \Pi(w)},
\end{align*}
where we defined $\rho_s := \frac{\underline{\lambda}(R)}{\bar{\lambda}(P)}$.

 Next, for the unstable modes, let $\sigma_v = a_i \in \Sigma_{v,a}$. 
Let $\tilde{\beta}>0$, and consider the auxiliary function $W_a = \tilde{\beta}x_e^\top x_e$. Recall that  $\tilde{\beta}(A_{a_i}^\top + A_{a_i}) = \tilde{\beta}\hat{R}_{a_i}$ as stated in Lemma \ref{lm:symR}. When $\vartheta \in \mathcal{D}$, during jumps we have that $W_a(x_e^+) = W_a(x_e)$. Moreover, when $\vartheta \in \mathcal{C}$, during flows we have:
\begin{align}\label{boundwunstable}
    \dot{W}_a(x_e) &= \tilde{\beta}(\dot{x}_e^\top x_e +  x_e^\top \dot{x}_e)\notag\\
    & = \tilde{\beta}\left(A_{a_i} x_e + A_{a_i}^{-1}B_{a_i}\dot{u} + A_{a_i}^{-1}E_{a_i}\Pi(w)\right)^\top x_e + \tilde{\beta}x_e^\top  \left(A_{a_i} x_e + A_{a_i}^{-1}B_{a_i}\dot{u} + A_{a_i}^{-1}E_{a_i}\Pi(w)\right)\notag\\
    & = \tilde{\beta}x_e^\top(A_{a_i}^\top  +  A_{a_i})x_e + 2 \tilde{\beta}x_e^\top (A_{a_i}^{-1}B_{a_i} \dot{u}) + 2 \tilde{\beta}x_e^\top (A_{a_i}^{-1}E_{a_i} \Pi(w))\notag\\
    & = \tilde{\beta}x_e^\top \hat{R}_{a_i} x_e - 2\varepsilon \tilde{\beta}x_e^\top  A_{a_i}^{-1}B_{a_i}(\varphi(u,w,x_e) + \varphi(u,w,0) - \varphi(u,w,0))+ 2 \tilde{\beta}x_e^\top (A_{a_i}^{-1}E_{a_i} \Pi(w))\notag\\
    &\leq \tilde{\beta}x_e^\top \hat{R}_{a_i} x_e + 2\varepsilon \tilde{\beta}\norm{x_e} \norm{ A_{a_i}^{-1}B_{a_i}}\norm{ \varphi(u,w,0) + (\varphi(u,w,x_e)  - \varphi(u,0))} + 2 \tilde{\beta}x_e^\top (A_{a_i}^{-1}E_{a_i} \Pi(w))\notag\\
    &\leq \tilde{\beta}x_e^\top \hat{R}_{a_i} x_e + 2\varepsilon \tilde{\beta}\norm{x_e} \norm{ A_{a_i}^{-1}B_{a_i}}\left(\norm{ \varphi(u,w,0)} + \norm{\varphi(u,w,x_e)  - \varphi(u,w,0)}\right)+ 2\tilde{\beta} x_e^\top (A_{a_i}^{-1}E_s \Pi(w))\notag\\
    & \leq \tilde{\beta}\bar{\lambda}(\hat{R}_{a_i})\norm{x_e}^2 + 2\varepsilon\tilde{\beta}\norm{ A_{a_i}^{-1}B_{a_i}}\norm{x_e}\norm{\nabla_u f(u,w)} + 2\varepsilon\ell_y\tilde{\beta}\norm{ A_{a_i}^{-1}B_{a_i}}\norm{C}\norm{G_{a_i}}\norm{x_e}^2\notag\\
    &\quad + 2 \tilde{\beta}\norm{x_e}\norm{  A_{a_i}^{-1}E_{a_i}}\norm{ \Pi(w)}\notag\\
    & \leq \rho_{a_i} W_a(x_e) +  2\varepsilon\tilde{\beta}\norm{ A_{a_i}^{-1}B_{a_i}}\norm{x_e}\norm{\nabla_u f(u,w)} + 2\varepsilon\ell_y\tilde{\beta}\norm{ A_{a_i}^{-1}B_{a_i}}\norm{C}\norm{G_{a_i}}\norm{x_e}^2\notag\\
    &\quad + 2\tilde{\beta} \norm{x_e}\norm{  A_{a_i}^{-1}E_{a_i}}\norm{ \Pi(w)},
\end{align}
where $\rho_{a_i} := \bar{\lambda}(\hat{R}_{a_i})$ for a fixed $a_i$. \blue{ We analyze the stability properties of the switched plant using the following auxiliary function:
\begin{equation}
\begin{aligned}
    \bar{W}(x_e,\tau_{v,1},\tau_{v,2},\sigma_v) = W_{\sigma_v}(x_e)e^\tau,~~\sigma_v\in \{s,a\},&  \quad \tau =  \ln(\omega)\tau_{v,1} + \tau_{v,2}(\rho_s + \rho_a),\\
    \text{where }0\leq \tau \leq \tau_0 := \ln(\omega)N_{0,v} + T_{0,v}(\rho_s + \rho_a)& \text{ and }~~\rho_s := \frac{\underline{\lambda}(R)}{\bar{\lambda}(P)}\qquad \rho_a :=  \bar{\lambda}(\hat{R}_a).
\end{aligned}
\end{equation}
By Assumption \ref{as:commonequilibrium}, $G_s =  G_{a_i}$, $\norm{ A_s^{-1}B_s}=\norm{ A_{a_i}^{-1}B_{a_i}}$, and  $\norm{ A_s^{-1}E_s}=\norm{ A_{a_i}^{-1}E_{a_i}}$ $~\forall a_i \in \Sigma_{v,a}$, hence, we drop the subindices $s$ and $a_i$, and we recall that $\bar{\lambda}(\hat{R}_a)$ was defined after Lemma \ref{lm:symR}. We set $\omega =\max\left\{\frac{\bar{\lambda}(P)}{\tilde{\beta}},\frac{\tilde{\beta}}{\underline{\lambda}(P)}\right\}$. 

The function $\bar{W}(x_e,\tau_{v,1},\tau_{v,2},\sigma_v) = W_{\sigma_v}e^\tau$ satisfies during jumps:
\begin{align}
\hspace{-0.5cm}
 \bar{W}(x_e^+,\tau_{v,1}^+,\tau_{v,2}^+,\sigma_v^+) = W_{\sigma_v^+}e^{\ln(\omega)(\tau_{v,1}-1)+\tau_{v,2}(\rho_s + \rho_a)} \leq \max_{\sigma_v^+ \in \Sigma_v}W_{\sigma_v^+}e^{-\ln(\omega)+\tau} \leq \omega W_{\sigma_v}e^{-\ln(\omega)}e^{\tau} = \bar{W}
\end{align}
} 
Similarly, during flows we have that:
\begin{align*}
\dot{\bar{W}}(x_e,\tau_{v,1},\tau_{v,2},\sigma_v) &= \dot{W}_s(x_e)e^\tau + W_s(x_e)e^\tau \dot{\tau}\\
    & = \dot{W}_s(x_e)e^\tau + W_s(x_e)e^\tau(\ln(\omega)\dot{\tau}_{v,1} + (\rho_s + \rho_a)\dot{\tau}_{v,2})\\
    & \leq \dot{W}_s(x_e)e^\tau +  W_s(x_e)e^\tau(\ln(\omega)\kappa_{v,1} + (\rho_s + \rho_a)\kappa_{v,2})\\
    & \leq -\rho_s W_s(x_e)e^\tau +  2\varepsilon\norm{PA^{-1}B}\norm{x_e}\norm{\nabla_u f(u,w)}e^\tau + 2\varepsilon\ell_ye^\tau\norm{PA^{-1}B}\norm{C}\norm{G}\norm{x_e}^2 \\
    & \quad + 2e^\tau \norm{x_e}\norm{ PA^{-1}E}\norm{ \Pi(w)} +  W_s(x_e)e^\tau(\ln(\omega)\kappa_{v,1} + (\rho_s + \rho_a)\kappa_{v,2})\\
    &\leq -\rho \bar{W}(x_e,\tau_{v,1},\tau_{v,2},\sigma_v) +  2\varepsilon\norm{PA^{-1}B}\norm{x_e}\norm{\nabla_u f(u,w)}e^\tau \\
    &\quad + 2\varepsilon\ell_ye^\tau\norm{PA^{-1}B}\norm{C}\norm{G}\norm{x_e}^2  + 2e^\tau \norm{x_e}\norm{ PA^{-1}E}\norm{ \Pi(w)}\\
&\leq -\rho \bar{W}(x_e,\tau_{v,1},\tau_{v,2},\sigma_v) +  2\varepsilon\ell\max\left\{\bar{\lambda}(P),\tilde{\beta}\right\}\norm{A^{-1}B}\norm{x_e}\norm{u - u^*_w}e^{\tau_0}\\
    &\quad + 2\varepsilon\ell_ye^{\tau_0}\max\left\{\bar{\lambda}(P),\tilde{\beta}\right\}\norm{A^{-1}B}\norm{C}\norm{G}\norm{x_e}^2 \\
    & \quad + 2e^{\tau_0} \max\left\{\bar{\lambda}(P),\tilde{\beta}\right\}\norm{A^{-1}E}\norm{x_e}\norm{ \Pi(w)},
\end{align*}
where $\rho: =  \rho_s - \kappa_{v,2}(\rho_s + \rho_a)-\ln(\omega)\kappa_{v,1}>0$, with $\omega \geq 1$. Similarly, using \eqref{boundwunstable}, we can obtain the following bound during unstable modes for the time-derivative of $\bar{W}(x_e,\tau_{v,1},\tau_{v,2},\sigma_v) =W_a e^\tau$:
\begin{align*}
\dot{\bar{W}}(x_e,\tau_{v,1},\tau_{v,2},\sigma_v) &= \dot{W}_a(x_e)e^\tau + W_a(x_e)e^\tau \dot{\tau}\\
& = \dot{W}_a(x_e)e^\tau + W_a(x_e)e^\tau(\ln(\omega)\dot{\tau}_{v,1} + (\rho_s + \rho_a)\dot{\tau}_{v,2})\\
& \leq \dot{W}_a(x_e)e^\tau +  W_a(x_e)e^\tau(\ln(\omega)\kappa_{v,1} +  (\rho_s + \rho_a)(\kappa_{v,2}-1))\\
& \leq \rho_{a} W_a(x_e)e^\tau +  2\varepsilon\tilde{\beta}e^\tau\norm{A^{-1}B}\norm{x_e}\norm{\nabla_u f(u,w)}+       2\varepsilon\ell_y\tilde{\beta}e^\tau\norm{ A^{-1}B}\norm{C}\norm{G}\norm{x_e}^2\\
    &\quad  + 2\tilde{\beta}e^\tau \norm{x_e}\norm{  A^{-1}E}\norm{ \Pi(w)} +  W_a(x_e)e^\tau(\ln(\omega)\kappa_{v,1} +  (\rho_s + \rho_a)(\kappa_{v,2}-1))\\
& \leq -\rho \bar{W}(x_e,\tau_{v,1},\tau_{v,2},\sigma_v) +  2\varepsilon\tilde{\beta}e^{\tau_0}\norm{A^{-1}B}\norm{x_e}\norm{\nabla_u f(u,w)}\\
&\quad + 2\varepsilon\ell_y\tilde{\beta}e^{\tau_0}\norm{A^{-1}B}\norm{C}\norm{G_{a_i}}\norm{x_e}^2 + 2\tilde{\beta}e^{\tau_0} \norm{x_e}\norm{  A_{a_i}^{-1}E_{a_i}}\norm{ \Pi(w)}\\
& \leq -\rho \bar{W}(x_e,\tau_{v,1},\tau_{v,2},\sigma_v) +  2\varepsilon \ell \max\left\{\bar{\lambda}(P),\tilde{\beta}\right\}\norm{A^{-1}B}\norm{x_e}\norm{u - u^*_w}e^{\tau_0}\\
    &\quad +2\varepsilon\ell_y e^{\tau_0}\max\left\{\bar{\lambda}(P),\tilde{\beta}\right\}\norm{A^{-1}B}\norm{C}\norm{G}\norm{x_e}^2\\
    &\quad  + 2e^{\tau_0} \max\left\{\bar{\lambda}(P),\tilde{\beta}\right\}\norm{  A^{-1}E}\norm{x_e}\norm{ \Pi(w)}.
\end{align*}
\noindent
To analyze the nominal controller \eqref{eq:gradientController} we considered the auxiliary function $V(u,w) = f(u,w) - f(u^*_w,w)$, then we get the following upper bound:
\begin{align*}
    \dot{V}(u,w) &= - \varepsilon \nabla_u f(u,w)^\top \varphi(x_e,u) + H^\top(\nabla f_y(Gu + Hw) - \nabla f_y(G u^* + Hw))\Pi(w)\\
    &= - \varepsilon \varphi(0,u)^\top(\varphi(0,u) + \varphi(x_e,u) - \varphi(0,u)) + \ell_y \norm{H}\norm{G}\norm{u - u^*}\norm{\Pi(w)}\\
    & \leq - \varepsilon\norm{\nabla_u f(u,w)}^2 + \varepsilon\norm{\nabla_u f(u,w)}\norm{G}\norm{\nabla f_y (Cx_e + Gu +Hw) - \nabla f_y (Gu +Hw)}\\
    & \quad + \ell_y \norm{H}\norm{G}\norm{u - u^*}\norm{\Pi(w)}\\
    & \leq - \varepsilon\norm{\nabla_u f(u,w)}^2 + \varepsilon\ell_y\norm{G}\norm{C}\norm{\nabla_u f(u,w)}\norm{x_e}+ \ell_y \norm{H}\norm{G}\norm{u - u^*}\norm{\Pi(w)}\\
    & \leq - \varepsilon\norm{\nabla_u f(u,w)}^2 + \varepsilon\ell_y\norm{G}\norm{C}\norm{\nabla_u f(u,w)}\norm{x_e}+ \ell_y \norm{H}\norm{G}\norm{u - u^*}\norm{\Pi(w)}\\
    & \leq - \varepsilon\mu^2\norm{u-u^*_w}^2 + \varepsilon\ell\ell_y\norm{G}\norm{C}\norm{u - u^*_w}\norm{x_e}+ \ell_y \norm{H}\norm{G}\norm{u - u^*}\norm{\Pi(w)}.
\end{align*}
Finally, for the HDS \eqref{eq:intercon_v}, we consider the Lyapunov function \blue{$U(\vartheta) =  \theta \bar{W}(x_e,\tau_{v,1},\tau_{v,2},\sigma_v) + (1-\theta) V(u,w) $, with $\theta \in (0,1)$. This function is bounded as follows:
\begin{align}\label{eq:quadbounds_Uv}
   \min\{\theta \underline{\lambda}(P),\theta \tilde{\beta},(1-\theta)\mu/2\}|\vartheta|^2_{\mathcal{A}^*}\leq U(\vartheta)\leq \max\left\{\theta\bar{\lambda}(P)e^{\tau_0},\theta\tilde{\beta}e^{\tau_0},(1-\theta)\ell^2/(2\mu)\right\}|\vartheta|^2_{\mathcal{A}^*}, 
\end{align}
$ \forall \vartheta \in \mathcal{C} \cup \mathcal{D}$, }where we used the facts that $\tau_{v,1} \in [0,N_{0,v}]$, $\tau_{v,2} \in [0,T_{0,v}]$, $\sigma_v \in \Sigma_v$, and $w \in \mathcal{W}$. For every $\vartheta \in \mathcal{D}$ we have that the jumps satisfy \blue{$U(\vartheta^+) \leq U(\vartheta)$}. Moreover, for $\vartheta \in \mathcal{C}$ the time derivative of $U$ is upper bounded as follows:
\begin{align*}
    \dot{U}(\vartheta) &= \theta\dot{\bar{W}}(x_e,u,\tau_{v,1},\tau_{v,2},\sigma_v) + (1-\theta) \dot{V}(u,w)\\
    &\leq -\rho \theta\min\{\underline{\lambda}(P),\tilde{\beta}\}\norm{x_e}^2 +  2\varepsilon \ell \theta\max\left\{\bar{\lambda}(P),\tilde{\beta}\right\}\norm{A^{-1}B}\norm{x_e}\norm{u - u^*_w}e^{\tau_0}\\
    &\quad +2\varepsilon\ell_y\theta e^{\tau_0}\max\left\{\bar{\lambda}(P),\tilde{\beta}\right\}\norm{A^{-1}B}\norm{C}\norm{G}\norm{x_e}^2 + 2\theta e^{\tau_0} \max\left\{\bar{\lambda}(P),\tilde{\beta}\right\}\norm{  A^{-1}E}\norm{x_e}\norm{ \Pi(w)}\\
    &\quad - (1-\theta)\varepsilon\mu^2\norm{u-u^*_w}^2 + \varepsilon\ell\ell_y(1-\theta)\norm{G}\norm{C}\norm{u - u^*_w}\norm{x_e}\\
    &\quad + \ell_y (1-\theta)\norm{H}\norm{G}\norm{u - u^*}\norm{\Pi(w)}.
\end{align*}
\noindent
Using $\xi : = (\norm{x_e}, \norm{u - u^*_w})$ we obtain
\begin{align}\label{eq:ISS_v1}
    \dot{U}(\vartheta) \leq -\varepsilon \xi^\top \Xi \xi + r^\top \xi \norm{\Pi(w)},
\end{align}
where $r: = \left(2\theta e^{\tau_0} \max\left\{\bar{\lambda}(P),\tilde{\beta}\right\}\norm{A^{-1}E}, \ell_y (1-\theta)\norm{H}\norm{G}\right)$, and $\Xi$ is a symmetric matrix as in Lemma \ref{lm:posdefmatrix} with parameters $\alpha:= \rho \min\{\underline{\lambda}(P),\tilde{\beta}\}$, $\beta:= 2\ell_y e^{\tau_0}\max\left\{\bar{\lambda}(P),\tilde{\beta}\right\}\norm{A^{-1}B}\norm{C}\norm{G}$, $\delta := \ell \ell_y\norm{C}\norm{G}$, $\chi:= 2\ell e^{\tau_0}\max\left\{\bar{\lambda}(P),\tilde{\beta}\right\}\norm{A^{-1}B}$, $\gamma:=  \mu^2$.  When $\theta:=\frac{\delta}{\delta + \chi}$, and $\varepsilon$ satisfies \eqref{eq:eps_v} we obtain that $\Xi$ is positive definite via Lemma \ref{lm:posdefmatrix}. Finally, we show that the the quadratic term dominates the linear term in \eqref{eq:ISS_v1}. Indeed, for $k \in (0,1)$, we can rewrite \eqref{eq:ISS_v1} as follows:
\begin{align}\label{eq:ISS_v2}
     \dot{U}(\vartheta) &\leq -\varepsilon(1-k)\underline{\lambda}(\Xi)\norm{\xi}^2 -\varepsilon k\underline{\lambda}(\Xi)\norm{\xi}^2+ r^\top \xi \norm{\Pi(w)}\nonumber\\
    &\leq -\varepsilon(1-k)\underline{\lambda}(\Xi)\norm{\xi}^2 -\varepsilon k\underline{\lambda}(\Xi)\norm{\xi}^2+ \norm{r}\norm{\xi} \norm{\Pi(w)}\nonumber\\
    &\leq -\varepsilon(1-k)\underline{\lambda}(\Xi)\norm{\xi}^2=-\varepsilon (1-k)\underline{\lambda}(\Xi) |\vartheta|_{\mathcal{A}^*}^2,~~~\forall~|\vartheta|_{\mathcal{A}^*} \geq \frac{\norm{r}}{\varepsilon \underline{\lambda}(\Xi) k }\norm{\Pi(w)}_t.
\end{align}
The result follows by applying \cite[Thm. 3.1]{CAI200947}, exploiting the quadratic bounds in \eqref{eq:quadbounds_Uv}, and the average dwell-time condition. 
\hfill $\blacksquare$


\section{Numerical Examples}
\label{section_numerical_examples}
In this section, we present two numerical experiments to illustrate our theoretical contributions. The first simulation utilizes a synthetic system; in the second experiment, we  consider  the problem of frequency regulation in power transmission systems. 
\subsection{Synthetic system}
We simulate a linear time-invariant system with state $x \in \R^2$, static control input $v \in \R$, dynamic control input $u \in \R$, and exogenous signal $w \in \R$. The goal is to regulate the solutions of the plant towards the solution of the following optimization problem
\begin{align}
   \min_{u,y} f_u(u) + f_y(y) := \min_{u,y} u^\top R u + (y - y_{\text{ref}})^\top Q (y - y_{\text{ref}}),
\end{align}
with $R = 2$, and $Q = [1~0;0~2]$.
\subsubsection{Attacks on the static feedback controller}
\begin{figure}
    \centering
    \includegraphics[width=0.45\columnwidth]{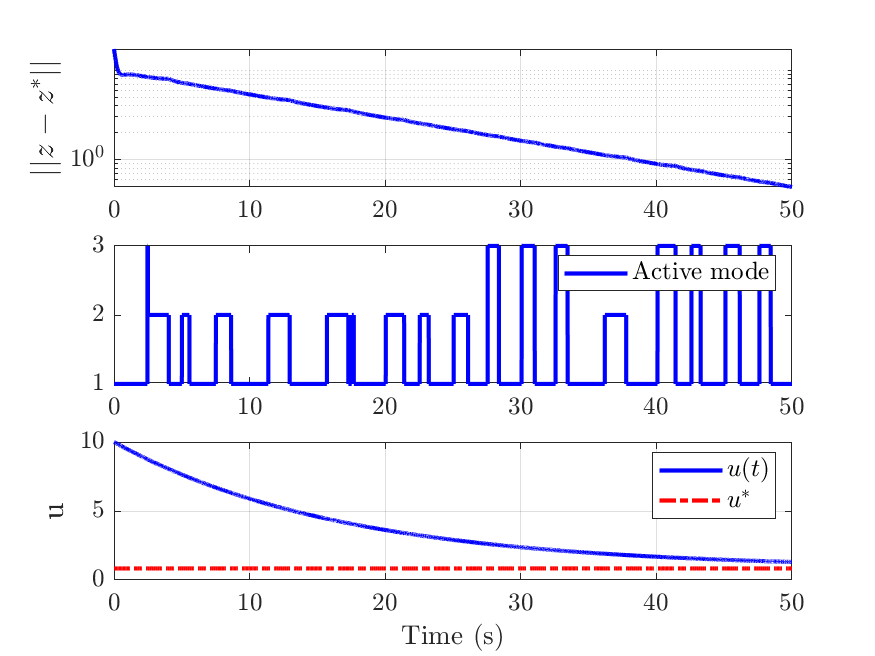}\includegraphics[width=0.45\columnwidth]{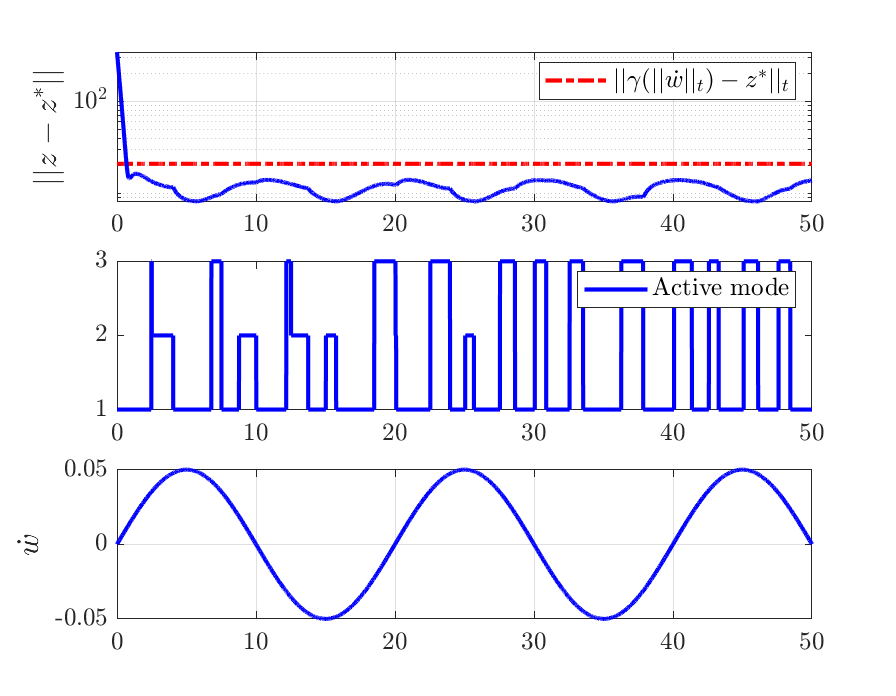}
    \caption{Attack to static feedback controller under constant and time-varying perturbations $w$.}
    \label{fig:static_1_wctt}
\end{figure}

We first consider the scenario where the static controller operates under two switching attacks, and the dynamic controller operates in the nominal stable mode. The plant has the following structure:
\begin{align}\label{eq:LTI_static1}
    \dot{x} &= \begin{bmatrix}1&0\\2&-1.5\end{bmatrix}x
    + \begin{bmatrix}1\\1\end{bmatrix}v
    +\begin{bmatrix}-1.06\\-0.62\end{bmatrix}u  +\begin{bmatrix}-0.82\\-0.79\end{bmatrix}w, \qquad
    y =  \begin{bmatrix}0.1&0\\0&0.1\end{bmatrix}x.
\end{align}
For this plant, we design an internal static controller that places the poles of the closed loop system $(F + NKC)$ at $-2$. Hence, $K = [-40 \quad 5]$. The attacks acting on the controller are modeled by the scalars $L_{a_1} = 0$, and $L_{a_2} = -0.1$, which leads to the following matrices:
\begin{align*}
    A_{a_1} &= \begin{bmatrix}1&0\\2&-1.5\end{bmatrix}\quad 
    B_{a_1} = \begin{bmatrix}0.34\\0.78\end{bmatrix}\quad  
    E_{a_1} = \begin{bmatrix}0.30\\0.34\end{bmatrix},\\
    A_{a_2} &= \begin{bmatrix}1.4&-0.5\\2.4&-1.55\end{bmatrix}\quad 
    B_{a_2} = \begin{bmatrix}0.48\\0.92\end{bmatrix}\quad  
    E_{a_2} = \begin{bmatrix}0.42\\0.45\end{bmatrix},
\end{align*}
which satisfy Assumption \ref{as:commonequilibrium}. The time-ratio parameter is set as $\kappa_{v,2} = 0.365$, which induces $\alpha = 1.71$ in Proposition \ref{pro:kappa_v2}. The gain inducing the time-scale separation is set as $\varepsilon^* = 0.0149$. In our simulation, the nominal stable mode is denoted by the index 1, the unstable mode $\sigma_{a_1}$ is denoted by the index 2, and $\sigma_{a_2}$ is denoted by the index 3. Thus $\Sigma_v=\{1,2,3\}$ and $\Sigma_{v,a}=\{2,3\}$. The theoretical value of $\varepsilon^*$ is conservative, and to obtain faster convergence we also simulated the system with a gain of  $\varepsilon = 20\varepsilon^*$. The left-hand side of Figure \ref{fig:static_1_wctt} shows the performance of the closed-loop system under a constant signal $w = 0.96$. As expected, the trajectories of the system converge to the set of optimal solutions. On the other hand, the right-hand side shows the performance of the system under a time-varying exogenous signal modeled by $\dot{w} = a \sin(\omega t)$, with $a = 0.05$ and $\omega = 2\pi \times 0.05$. Here, we can see that the norm of the tracking error is eventually bounded by 21. 
\begin{figure}
    \centering
    \includegraphics[width=0.45\columnwidth]{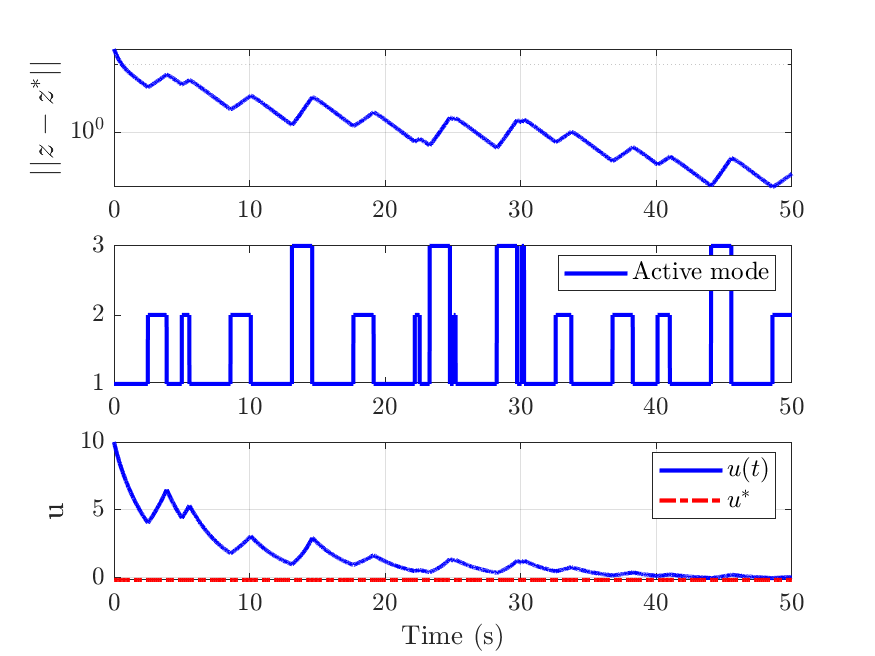}\includegraphics[width=0.45\columnwidth]{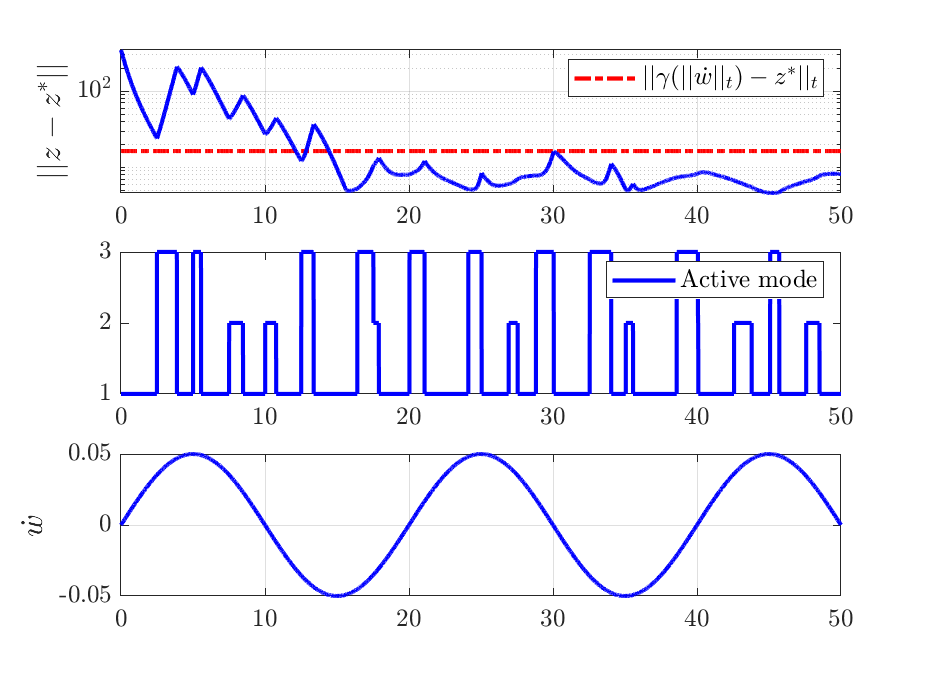}
    \caption{Attack to dynamic feedback controller under constant and time-varying perturbations $w$.}
    \label{fig:dynamic_1_wctt}
\end{figure}
\subsubsection{Attacks on the dynamic feedback controller}
We now consider the dynamic controller operating under two switching attacks, while the nominal static controller (stable) remains free from any attack. In this case, we have the following nominal plant
\begin{align}\label{eq:LTI_dyn1}
    \dot{x} &= \begin{bmatrix}-3&0.5\\-2&-1\end{bmatrix}x
    +\begin{bmatrix}1\\1\end{bmatrix}u  
    +\begin{bmatrix}1\\1\end{bmatrix}w, \qquad
    y =  \begin{bmatrix}1&0\\0&1\end{bmatrix}x.
\end{align}
We simulate the system in \eqref{eq:LTI_dyn1} interconnected with the dynamic controller \eqref{eq:attackModelDynamic} under attacks. The attacks are modeled by the scalar gains $M_{a_1} = -1$, and $M_{a_2} = -2$. The time-ratio parameter is set to $\kappa_{u,2} = 0.33$. To induce time-scale separation, the gain is set as $\varepsilon^* = 0.0093$. In the simulation, the nominal stable mode is characterized by the index $\{1\}$, the attack $\sigma_{a_1}$ is characterized by the mode $\{2\}$, and also $\sigma_{a_2}$ is characterized by \{3\}. Since simulations show that the theoretical gain $\varepsilon^*$ is very conservative, we set $\varepsilon = 20 \varepsilon^*$ to achieve faster convergence in the experiments. In the left plot of Figure \ref{fig:dynamic_1_wctt}, we show the decreasing tracking error under a constant $w$. Similarly, in the right plot, we show the behavior of the system under a time-varying signal $\dot{w} = a\sin(\omega t)$ with $a = 0.05$ and $\omega  = 2\pi \times 0.05$. As predicted by Theorem \ref{thm:iss_v}, the tracking error is eventually uniformly ultimately bounded.

\subsection{Persistent Attack in a Power System}
We study the problem of regulating frequency in a power transmission system; in particular, we model the frequency dynamics as an LTI system based on a DC model of the network and a linearized model of the swing equation; see, e.g., the experimental setup in \cite[Section V]{MC-ED-AB:19}.\footnote{In this section, we denote as $I_n$  the identity matrix of dimension $n$, the vector of ones as $\mathbb{1}$, $0_{n \times m}$ is a matrix of zeros of dimension $n \times m$, and the symbol diag() transforms a vector to a matrix with the elements of the vector on the diagonal.} Accordingly, we consider the following LTI dynamical system:
\begin{equation}\label{eq:swingeq}
\begin{aligned}
    \begin{bmatrix}\dot{\mathbf{\theta}}\\\dot{\mathbf{\omega}}\end{bmatrix} &= \underbrace{\begin{bmatrix}0 & I \\ -M^{-1} Y & -M^{-1}D\end{bmatrix}}_{\bar{A}}\begin{bmatrix}\mathbf{\theta}\\\mathbf{\omega}\end{bmatrix} + \underbrace{\begin{bmatrix}0\\B_1\end{bmatrix}}_{\bar{B}}u + \underbrace{\begin{bmatrix}0\\E_1\end{bmatrix}}_{\bar{E}}w,\\
    \underbrace{\begin{bmatrix}p_l\\\bar{\omega}\end{bmatrix}}_y &= \underbrace{\begin{bmatrix}Y & 0 \\ 0 & \frac{1}{n}\mathbb{1}^\top\end{bmatrix}}_{\bar{C}}\begin{bmatrix}\mathbf{\theta}\\\mathbf{\omega}\end{bmatrix},
\end{aligned}
\end{equation}
where the state $z := (\theta,~\omega) \in \mathbb{R}^{2n_1}$ collects the phase angles $\theta_i$ and frequencies $\omega_i$ at $n_1$ nodes, the vector $u$ collects the powers generated by generator (modeled as synchronous machines),  $w$ represents uncontrollable loads, $p_l$ is the vector of power flows on the lines, and $\bar{\omega}$ is the average system frequency. The matrix $Y \in \mathbb{R}^{n_1 \times n_1}$ is the susceptance matrix of the grid whose row sum is equal to zero, i.e. $Y\mathbb{1} = 0$. The matrix $B_1 \in \{0,1\}^{n_1 \times m}$ sets the controllable nodes or generators i.e. $B_1(i,j) = 1$ when $i$ is a generator; and $E_1 \in \{0,1\}^{n_1 \times q}$ sets the uncontrollable loads, i.e. $B_2(i,j) = 1$ when $i$ is a load. Finally,  $M =  \text{diag}(m_1,\ldots,m_{n_1})$ and $D =  \text{diag}(d_1,\ldots,d_{n_1})$ are diagonal matrices respectively containing the rotational inertia coefficients and the damping coefficients in the diagonal. 

System \eqref{eq:swingeq} is marginally stable since $\bar{A}$ has an eigenvalue at zero with eigenvector $(\mathbb{1}, 0)$. We eliminate the average mode $\bar{\theta}$ with the following coordinate transformation~\cite[Section V]{MC-ED-AB:19}:
\begin{align}
    x := \begin{bmatrix}\theta \\ \omega\end{bmatrix} = \underbrace{\begin{bmatrix}U & 0 \\ 0 & I\end{bmatrix}}_{T} x + \begin{bmatrix}\mathbb{1}\\ 0\end{bmatrix}\bar{\theta},
\end{align}
where the columns of the matrix $U \in \R^{n_1 \times (n_1 - 1)}$ form an orthonormal basis orthogonal to span($\mathbb{1}$). One way to select this $U$ is taking the $n_1 -1$ eigenvectors of $Y$ that correspond to the nonzero eigenvalues. The new system with state $x =  T^\top x \in \R^{2n_1 -1}$ takes the following form
\begin{equation}\label{eq:LTI_ex2}
 \begin{aligned}
    \dot{x} &= Ax + Bu + Ew,\\
    y &= Cx,
\end{aligned}  
\end{equation}
where $A := T^\top \bar{A}T$, $B:= T^\top \bar{B}$, $E:= T^\top \bar{E}$, and $C:=\bar{C}T$. System \eqref{eq:LTI_ex2} is stable and the output $y$ remains unaltered by the coordinate transformation. Based on this, we consider the following optimization problem associated with the transmission systems:
\begin{equation}\label{eq_optpro_ex2}
\begin{aligned}
    \min_{u} & \qquad u^\top Q u + f_y(y)\\
    \text{s.t}& \qquad y = Gu + Hw
\end{aligned}
\end{equation}
where $G = -C A^{-1}B$, $H = -CA^{-1}E$, $Q \in \R^{m \times m}$ is positive definite and $f_y(y) = \eta \max \left\{\underline{y}_i - y_i, y_i - \overline{y}_i,0\right\}^2$. The function $f_y(\cdot)$ imposes soft constraints to the system output $y$; moreover, $f_y(y)$ is continuously differentiable with a  $\eta-$Lipschitz gradient. With this formulation, the dynamic feedback controller is given by: 
\begin{align}\label{eq:ctrl_ex2}
    \dot{u} = -\varepsilon(2Hu + \eta G^\top \nabla f_y(y)),
\end{align}
where $\varepsilon >  0$, and $\nabla f_y(\cdot)$ boils down to the following soft thresholding operator:
\begin{align*}
    \nabla f_y(v) = \begin{cases}
            v_i -  \underline{y}_i, & v_i \leq \underline{y}_i \\
            0 & \underline{y}_i < v_i < \overline{y}_i\\
            v_i - \overline{y}_i & v_i \geq \overline{y}_i. 
    \end{cases}
\end{align*}
The nominal system that we simulate contains a simplified version of the IEEE 9Bus Test case; see \cite{MC-ED-AB:19} for a schematic of the network. In particular, it features  $n_1 = 6$ nodes, with generators at $3$ nodes, and uncontrollable loads under disturbances in the remaining $3$ (with loads modeled as ``negative'' generators). We build the matrices in system \eqref{eq:swingeq} as follows  $M = \text{diag}(1,1,1,0.1,0.1,0.1)$, $D = 0.1*I_{n_1}$, $B_1 = \text{diag}(1,1,1,0,0,0)$, and $E_1 = (0_{(n_1 - q) \times q},I_{q})$, and 
\begin{align*}
    Y =  \begin{bmatrix}
        13.5776     &    0     &    0  & -6.9783  & -6.5993     &    0\\
          0  & 11.7898    &     0  & -4.3844   &      0  & -7.4054\\
          0     &    0 &  10.4895     &    0 &  -4.2507  & -6.2388\\
    -6.9783  & -4.3844     &    0  & 11.3627     &    0    &     0\\
    -6.5993     &    0  & -4.2507     &    0 &  10.8500     &    0\\
          0 &  -7.4054  & -6.2388     &    0      &   0  & 13.6441
    \end{bmatrix}.
\end{align*}
Since $A \in \R^{n \times n}$, where $n := 2n_1 - 1$, we use $R = I_n$ to compute $P$ according to Assumption \ref{as:stabilityPlant}. We consider a constant $w =(-0.9, -1 , -1.25)$ for the first $100~s$, and $-0.45(\cos(2\pi 0.05 t ) +1)$ for the remaining $100~s$. For the cost function, we set $Q =  I_m$; furthermore, the parameters for $f_y$ are $\eta =10$, $\overline{y} = 5*\mathbb{1}$, and $\underline{y} = -5*\mathbb{1}$.

\begin{figure}
    \centering
    \includegraphics[width=0.45\columnwidth]{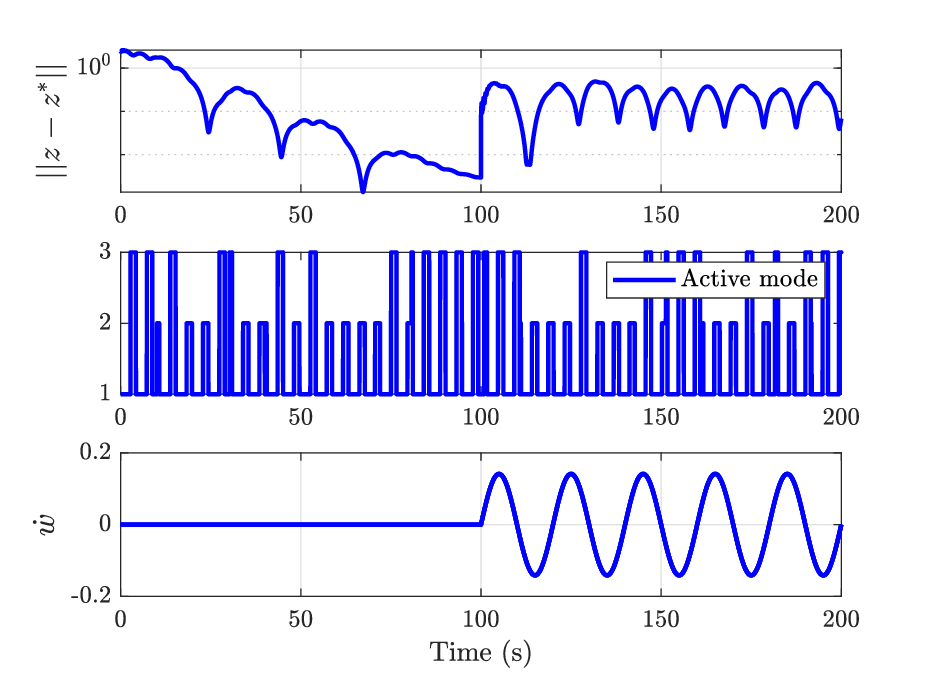}
    \includegraphics[width=0.45\columnwidth]{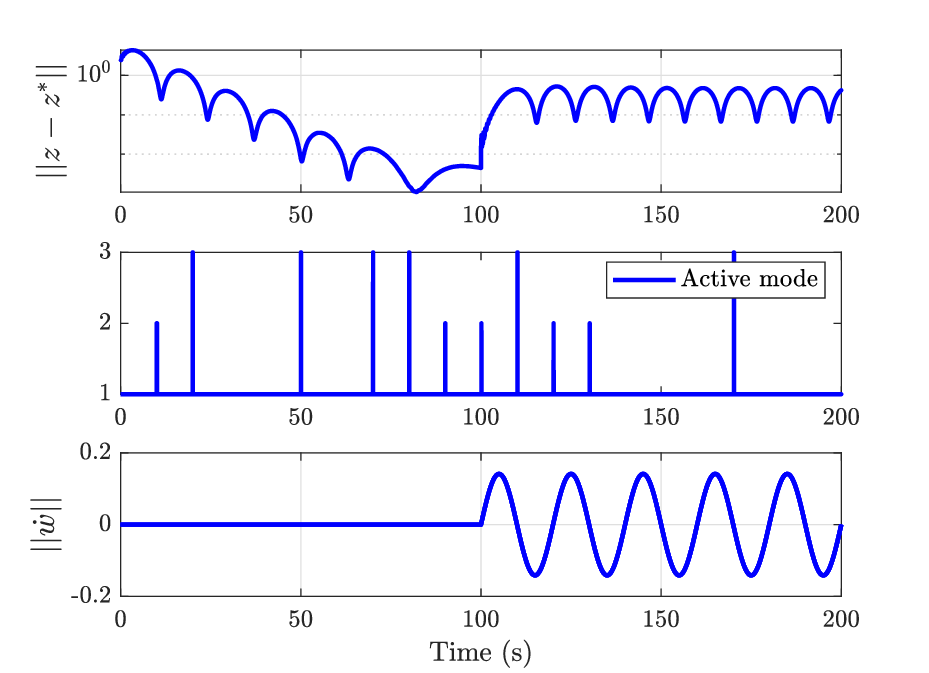}
    \caption{
    In the right figure, we present an attack to dynamic feedback controller to system in \eqref{eq:LTI_ex2} with $\varepsilon = 1e5\varepsilon^*$. In the left figure, we present an attack to the plant in \eqref{eq:LTI_ex2}.
    }
    \label{fig:ex2}
\end{figure}

We simulate attacks to the dynamical controllers with $M_{a_2} = -1$, and $M_{a_3} = -2$. 
The right plots of Figure \ref{fig:ex2} show the results of this simulation. We implement our algorithm \eqref{eq:ctrl_ex2} with $\varepsilon =  1e5\varepsilon^*$; we note that the error $||z - z^*||$ decreases whenever $\norm{\dot{w}}$ is constant, and remains bounded whenever $w$ is time-varying. Moreover, we simulate attacks acting on the plant by systematically modifying the matrix $M$ in the following way: $M_2 = \text{diag}(1,1,-101,0.1,0.1,0.1)$, and $M_3 = \text{diag}(1,1,-105,0.1,0.1,0.1)$. This modifications lead to an unstable matrix $A$. On the left side of Figure \ref{fig:ex2} we present the results. As expected, the error $||z - z^*||$ decreases when $w$ remains constant, but only remains bounded whenever $w$ is time-varying. These results are consistent with our theoretical results.

\section{Conclusions}
\label{section_conclusions}
In this paper, we presented the first systematic stability analysis of LTI plants controlled by dynamic gradient-based controllers in subject to persistent attacks. For cost functions that satisfy certain smoothness and gradient-dominated properties, we showed that the static output-feedback controllers guarantee exponential stability of the closed-loop system whenever the total activation time of the attacks in any given window of time is sufficiently bounded, and not too frequent. For the dynamic gradient-flow controllers, we showed that system stability can be further guaranteed by properly tuning the controller gain. To the best knowledge of the authors, our results are the first in the literature of feedback-based automatic optimization that study the convergence and stability properties of the algorithms using the framework of hybrid dynamical systems, and Lyapunov-based tools for switched systems with unstable modes. Overall, our results demonstrate that, under exogenous disturbances, input-to-state stability can be guaranteed for the closed-loop system provided the total activation time of the unstable modes satisfies a particular upper bound, and enough time-scale separation is induced between the dynamics of the plant and the controller. We also presented explicit characterizations of the ISS gains, the theoretical upper bounds on the controller gains, and the time ratio constraints needed to preserve stability. Our results were validated in two numerical examples. Relevant research directions that require further investigation include the analysis of switching plants with non-common equilibrium points and extensions to plants with nonlinear dynamics.

\bibliographystyle{model1-num-names}

\bibliography{alias.bib,biblio.bib}

\appendix
\vspace{-0.1cm}
\section*{Appendix}

\begin{lemma}\label{lm:posdefmatrix}
Let $\alpha, \beta, \delta, \chi, \gamma$ be positive scalars, let $\theta \in (0,1)$ be tunable parameters, and let
\begin{align*}
    \Xi = \begin{bmatrix}
        \theta\left(\frac{\alpha}{\varepsilon} - \beta\right) & -\frac{1}{2}\left((1-\theta)\delta + \theta \chi\right)\\
        -\frac{1}{2}\left((1-\theta)\delta + \theta \chi\right) & (1-\theta)\gamma
    \end{bmatrix}.
\end{align*}
If $0 < \varepsilon < \alpha \gamma / (\beta \gamma + \delta \chi)$, then there exists $\theta < \delta / (\delta + \chi)$, such that $\Xi$ is positive definite.
\end{lemma}
\textbf{Proof:} The matrix $\Xi$ is positive definite if and only if the leading principal minors are positive. In this case $(1 - \theta)\gamma > 0$ and $\theta(1-\theta)\left(\frac{\alpha}{\varepsilon} - \beta\right)\gamma > \frac{1}{4}((1-\theta)\delta + \theta \chi)^2$. The first inequality is guarantee by the definition of $\theta$ and $\gamma$. The second inequality can be rewritten as:
\begin{align*}
    \varepsilon < \frac{\alpha \gamma }{\beta \gamma +\frac{((1-\theta)\delta + \theta \chi)^2}{4\theta(1-\theta)}} = \hat{\varepsilon}(\theta).
\end{align*}
The function $\hat{\varepsilon}$ attains its maximum at $\theta = \theta^* := \frac{\delta}{\delta + \chi}$, with $\varepsilon^* := \hat{\varepsilon}(\theta^*) = \frac{\alpha \gamma}{\beta \gamma +  \delta \chi}$. \hfill $\blacksquare$

\end{document}